\documentstyle{jtbbook}

\makeatletter
\def\@begintheorem#1#2{\it \trivlist \item[\hskip
\labelsep{\bf #2\ #1.}]}
\def\@opargbegintheorem#1#2#3{\it \trivlist
      \item[\hskip \labelsep{\bf #1\ #2\ (#3).}]}
\makeatother

\newtheorem{them}{Theorem}[section]
\newtheorem{quest}{Question}
\newtheorem{lemm}[them]{Lemma}
\newtheorem{cor}[them]{Corollary}
\newtheorem{prop}[them]{Proposition}
\newtheorem{thm}[them]{Theorem}

\newcommand{\jtbnumpar}[1]{\refstepcounter{them}
\trivlist
\item[\hskip \labelsep{\bf \thethm \ #1.}]}

\newcommand{\jtbdef}{\jtbnumpar{Definition}}
\def\jtbnot{\jtbnumpar{Notation}}
\let\endjtbdef=\par
\let\endjtbnumpar=\par
\let\endjtbnot=\par

\newtheorem{ssnote}{COMMENT}

\def\PROOF #1.{\par\noindent{\it Proof#1}.\ \ignorespaces}
\def\proof #1.{\par\noindent{\it Proof#1}.\ \ignorespaces}

\def\subm{\leq}
\def\extm{\geq}

\def\Mscr{{\cal M}}

\let\?=\joinrel

\let\leftv=^
\let\rightv=^

\def\endproof{\vskip 5pt plus1pt minus1pt }

\def\emptynode{{\scriptstyle\langle\;\rangle}}

\def\k/{\kern.2em}    

\def\NF{\mathop{\rm NF}}

\def\LFP{\mathop{\rm LFP}}

\def\cpr{\mathop{\rm cpr}}

\def\im{\mathop{\rm rng}}

\def\ls{\rm LS}               %

\let\bet=\Im        %

 \mathchardef\bet="0B69

\def\down{\smash{\mathchar"0223}}

\let\union=\cup             %

\edef\bigcup{\mathop{\textstyle\mathchar\the\bigcup}}

\let\inter=\cap             %

\edef\bigcap{\mathop{\textstyle\mathchar\the\bigcap}}

\edef\bigwedge{\mathop{\textstyle\mathchar\the\bigwedge}}

\edef\bigvee{\mathop{\textstyle\mathchar\the\bigvee}}

\edef\sum{\mathop{\textstyle\mathchar\the\sum}}
\def\ind #1#2#3{#1 \mathbin{\mathop{\down}_{#2}} #3}
\def\math&{\ \& \ }

\def\force {\mathrel^\joinrel\rightarrow}
\def\force {\mathrel{\scriptstyle\mathrel^\joinrel\rightarrow}}
\def\forceq {\mathrel{\mathop{\force}\limits_{\textstyle\texsim}}}
\def\forceq{\mathrel^\joinrel
 \mathrel{\mathop{\rightarrow}\limits_{\smash{\textstyle\texsim}}}}
\def\forceq{\mathrel{\scriptstyle\mathrel^\joinrel
 \mathrel{\smash{\mathop{\rightarrow}\limits_{\smash{\raise
 2pt\hbox{$\scriptstyle\texsim$}}}}}}}

\let\exclaim=!                 %
\let\texsim=\sim         %
\let\conj\sim             %
\def\conjp #1 {\conj_{#1}}     %
\let\sim\simeq             %
\let\neg=\lnot             %

\def\0bar{\bar 0}         %
\def\1bar{\bar 1}         %

\def\rom #1{\hbox{#1}}

\def\hbar{\overline h}

\def\sbar{\overline s}
\def\tbar{\overline t}
\def\ubar{\overline u}
\def\vbar{\overline v}

\def\Mun{\underline M}
\def\Nun{\underline N}

\def\emptynode{{\langle\;\rangle}}

\let\ocirc\circ    

\def\bet{]}
\def\ord{\rom ord}

\def\bet{]}

\def\Mun{{\underline M}}
\def\Nun{{\underline N}}

\def\text#1{\ifmmode\leavevmode\hbox{#1}\else
   \typeout{Warning: \string\text \space used outside math mode!}
   \begingroup\hbox{#1}\endgroup\fi}

\def\ord{\rom ord}
\title{The Primal Framework I.}
\author{J. Baldwin
\thanks{Partially supported by N.S.F. grant 8602558 }\\
Department of Mathematics\\
University of Illinois, Chicago
\and
S. Shelah
Department of Mathematics\\
Hebrew University of Jerusalem
\thanks{Both authors thank Rutgers University and the U.S. Israel
Binational Science foundation for their support of this project.
This is item 330 in Shelah's bibliography.}}
\date{December 17, 1987; revised February 21, 1989}
\begin{document}
\maketitle
This the first of a series of articles dealing with abstract
classification theory.  The apparatus to assign systems of cardinal
invariants to models of a first order theory (or determine its
impossibility) is developed in \cite{Shelahbook}.  It is natural to try
to extend this theory to classes of models which are described in other
ways.  Work on the classification theory for nonelementary classes
\cite{Shelahnonelemii} and
for universal classes \cite{Shelahuniversal} led to the conclusion
that an axiomatic approach provided the best setting for developing a
theory of wider application.
This approach is reminiscent of the early work of Fraiss\'e and J\'onson
on the existence of homogeneous-universal models.
As this will be a long project it seems
appropriate to report our progress as we go along.

In large part this series of articles will parallel the development in
\cite{Shelahuniversal}.  A survey of that paper which could serve as
an introduction to this one
is \cite{Baldwinclass85}.
The first chapter
of this article corresponds to Section 2
of \cite{Shelahuniversal}.  In it we describe
the axioms on which the remainder of the
article depends and give some examples and context to justify this level
of generality.  As is detailed later the principal goal of this series
is indicated by its title.  The study of universal classes takes as a
primitive the notion of closing a subset under functions to obtain a
model.  We replace that concept by the notion of a prime model.  We begin
the detailed discussion of this idea in Chapter II.  One of the
important contributions of classification theory is the recognition that
large models can often be analyzed by means of a family of small models
indexed by a tree of height at most $\omega$.  More precisely, the
analyzed model is prime over such a tree.  Chapter III provides
sufficient conditions for prime models over such trees to exist.  The
discussion of properties of a class which guarantee that each model in
the class is prime over such a tree will appear later in the series.

We introduce in Chapters I and II a number of principles, which we
loosely refer to as axioms.  At the beginning of Chapter III we define
the notion of an adequate class --- a class which satisfies those axioms
that we assume in the mainline of the study.  This notion of an adequate
class will be embellished by further axioms in later papers of this
series.
Our use of the word axiom in this context is somewhat inexact; postulate
might be better.  In exploring an unknown area we list certain
principles which appear  to make important distinctions.  In our
definition of an adequate class we collect a family of these principles
that is sufficient establish a coherent collection of results.
We thank N. Shi for carefully reading the paper and making a number of
helpful suggestions.
\chapter{THE ABSTRACT FRAMEWORK}

Shelah developed in \cite{Shelahuniversal} several frameworks
for studying aspects
of classification theory.  In each case he studied a triple
$${\bf K} =
\langle K, \subm_{\bf K},{\NF}_{\bf K} \rangle;$$
$K$ is a collection of structures, $\subm_{\bf K}$ denotes elementary
submodel
with respect to {\bf K}, and $\NF_{\bf K}$ is a $4$-ary relation
(nonforking)
denoting that certain models are in stable amalgamation.  The
original paper primarily studied classes which admitted a fourth basic
notion:
`generated submodel'.  We generalize that context here by taking
as a fourth basic component a
predicate $\cpr$. Intuitively,
$\cpr =
\cpr_{\bf K}$
holds
of a structure $M\in {\bf K} $ and
a  chain of models $\Mun$ (Section I.2) if $M$ is prime over
$\Mun$. Thus this paper studies quadruples $${\bf K} =
\langle K, \subm_{\bf K},{\NF}_{\bf K}, {\cpr}_{\bf K}\rangle.$$

Section 1 reviews the properties of
elementary submodel and free amalgamation which carry over from
\cite{Shelahuniversal}.
In section 2 we provide a number of examples of classes which satisfy
the basic axioms.

\section{Basic properties of $\subm_{\bf K}$ and NF}

$K$ always denotes a class of structures of a fixed similarity type.
$K$ and all relations that we define on it are assumed to be closed
under isomorphism.
Although technically both $\subm$ and NF should be subscripted with {\bf
K},
we usually omit the subscript for ease of reading.
$M$ and  $N$  (with subscripts) denote members of $K$ unless we
explicitly say otherwise.  $A$ and $B$ will denote subsets of members of
$K$.
We write $M$ is contained in $N$ ($M \subseteq N$) if $M$ is a
substructure of $N$; that is,  if the universe of $M$ is a subset of
that of $N$, the relations of $N$ are those imposed by $M$ and
$M$ is closed under any functions in the language.
We write $M$ is a submodel or more explicitly a $K$-submodel of $N$ for
$M \subm N$.

The first group of axioms describe our notion of elementary submodel.
\jtbnumpar{Axiom Group A: {\bf K}-submodels}
\begin{description}
    \item
    \item[A0]  If $M \in K$ then $M \subm M$.
    \item[A1]  If $M \subm N$ then $M$ is a substructure of $N$
    \item[A2]  $\subm$ is transitive.
    \item[A3]  If $M_0 \subseteq M_1 \subseteq N$, $M_0 \subm N$ and $M_1
\subm N$ then $M_0 \subm M_1$.
\end{description}
\medskip
It is sometimes essential to distinguish between $M \subm_{\bf K} N$
which implies $M \subseteq N$ and the existence of an embedding $f$
of $M$
into $N$ whose image is a {\bf K}-submodel of $N$.
\jtbdef A {\em {\bf K}-embedding} is an isomorphism $f$ from an $M$ in
$K$ to an $N$ in $K$ such that $\im f \subm_{\bf K} N$.
$N$ is then a {\em {\bf K}-extension } of $M$ via $f$.

If $f$ is not mentioned explicitly then it is the identity.

We require one further important property of {\bf K} and $\leq_{\bf K}$.
\jtbdef
{\bf K} has the $\lambda$-L\"owenheim-Skolem
property
($\lambda$-LSP)
or (the
L\"owenheim-Skolem property down to $\lambda$) if
$A \subseteq M \in {\bf K} $ and $|A| \leq \lambda$
implies there is an $N \in {\bf
K} $ with $A \subseteq N \subm M$ and $|N|\leq \lambda$.
The
L\"owenheim-Skolem
number of {\bf K} is the least $\lambda$
such that {\bf K} has the $\lambda$-L\"owenheim-Skolem
property.
We write $\ls
({\bf K} )= \lambda$.

Note that the set of cardinals for which {\bf K} has the L\"owenheim
Skolem property may not be convex.  Moreover the related requirement on
$\lambda$, for any $A$ there is an $N\supseteq A$ with $N\in {\bf K} $
and $|N| \leq |A| + \lambda$, is still different and will be
investigated later.

\endjtbdef
\jtbnumpar{Axiom Group A: {\bf K}-submodels}
\begin{description}
    \item
    \item[A4]  $\ls ({\bf K} ) < \infty$.
\end{description}
\endjtbnumpar
Our axioms differ from those of Fraiss\'e and J\'onsson in that there is
neither a joint embedding nor an amalgamation axiom.  Our approach here
is to assume
in the next set of axioms a particularly strong form of amalgamation.
The use of nonamalgamation as a
source of nonstructure has been explored by Shelah in several places.
See especially Chapter I of \cite{Shelahuniversal} and
\cite{Makowskyabstractembed} and its progenitor \cite{Shelahnonelemii}.
We will obtain joint embedding by fiat (by restricting to a subclass
that satisfies it).  The J\'onsson Fraiss\'e constructions also require
closure under unions of chains.  This requirement is more subtle than it
first appears; it is the major topic of \cite{BaldwinShelahprimalii}.

We say two structures are {\em compatible} if they (isomorphic copies of
them) have a common {\bf
K}-extension. Any class
{\bf K}
that satisfies Axiom C2 (below)
is split into classes with the
joint
embedding property by the equivalence relation of compatibility and if
{\bf K} has a L\"owenheim-Skolem number then there will be only  a set
of equivalence classes.  Given any fixed model $M$ (or diagram of
models) this equivalence relation will be refined by `compatibility over
$M$'.  We explore this refinement in Chapter II.

The second group of axioms concern the independence relation.  We begin
by describing a relation of four members of $K$.
$$\NF(M_0,M_1,M_2,M_3)$$
means that $M_1$ and $M_2$ are freely amalgamated over $M_0$ within
$M_3$.  The notation $\NF$ arises from the reading the type of $M_1$
over $M_2$ inside $M_3$ does not fork over $M_1$.
We will eventually show a dichotomy between nonstructure results
and
the existence of a monster model $\Mscr$.  Thus, in trying to establish
a structure theory
we can introduce a $3$-ary relation $\NF(M_0,M_1,M_2)$ to
abbreviate $\NF(M_0,M_1,M_2,\Mscr)$.  We usually write this relation as
$\ind {M_1} {M_0} {M_2}$.  Even before showing the existence of the
monster model we will employ this notation if the choice of $M_3$ is
either clear from context or there are several possibilities which serve
equally well.

\jtbnot
A $4$-tuple $\langle M_0,M_1,M_2,M_3\rangle$
is called a {\em full free amalgam}
if it satisfies
$\NF(M_0,M_1,M_2,M_3)$.
The three-tuple
$\langle M_0,M_1,M_2\rangle$ is called a {\em free
amalgam} if it is an initial segment of a full free amalgam.
We often write $M_1$ and $M_2$ are freely amalgamated over $M_0$ in
$M_3$.  We refer to such a diagram as a `free vee'.
An isomorphism between two free amalgams $\Mun = \langle
M_0,M_1,M_2\rangle$
and $\Mun'=\langle M'_0,M'_1,M'_2\rangle$ is a triple of
isomorphisms $f_i$ mapping $M_i$ to $M'_i$ with $f_0$ contained in $f_1$
and $f_2$.
There is no guarantee (until we assume Axiom C5 below) that the
isomorphisms $f_1$ and $f_2$ have a common extension to an $M_3$ which
completes $\Mun$.
We extend this notion of isomorphism to arbitrary diagrams in Section 2.
\jtbnumpar{Axiom Group C: Independence}
The following axioms describe the independence relation.
For convenience of comparison we have kept the numbering from
\cite{Shelahuniversal} when we have just copied an axiom.  Some
axioms from that
list (e.g. C4) are omitted here.
In particular, the role of Axiom Group B from \cite{Shelahuniversal},
which dealt with the notion of
generation, is played here by Axiom Group D.
(See Section~\ref{primemodelsovervees}.)

\begin{description}
     \item[C1]  If $NF(M_0,M_1,M_2,M_3)$ then
$M_0 \subm M_2 \subm M_3$
and
$M_0 \subm M_1 \subm M_3$.  In particular, each $M_i \in K$.
If $\NF(M_0,M_1,M_2,M_3)$, we say $M_1$ and $M_2$ are freely amalgamated
(or independent) over $M_0$ in $M_3$.
     \item[C2 Existence]  If $M_0$ is a $K$-submodel of both $M_1$ and
$M_2$ then there are copies
(over $M_0$)
$M'_1$ and $M'_2$
of
$M_1$ and $M_2$ which are freely amalgamated in some $M_3 \in K$.
     \item[C3 Monotonicity]
        \begin{enumerate}
           \item  If $M_1$ and $M_2$ are freely amalgamated over $M_0$
in $M_3$ then so are $M_1$ and $M'_2$ for any $M'_2$ with $M_0 \subm
M'_2 \subm M_2$.
           \item  If $M_1$ and $M_2$ are freely amalgamated over $M_0$
in $M_3$ then they are freely amalgamated in any $M'_3 \extm M_3$.
           \item  If $M_1$ and $M_2$ are freely amalgamated over $M_0$
in $M_3$ then they are freely amalgamated in any $M'_3$ containing
$M_1 \union M_2$ and with $M_3' \leq M_3$.
\end{enumerate}
     \item[C5 Weak Uniqueness]
Suppose
$\langle \Mun, M_3\rangle$ and
$\langle \Mun', M'_3\rangle$
are full free amalgams.
If $\Mun$ and $\Mun'$ are isomorphic free amalgams (via $f$)
then there is an $N\in {\bf K} $ with
$M_3'\subm N$,
and an extension of $f$ mapping $M_3$ isomorphically
onto a {\bf K}-submodel of $N$.
     \item[C6 Symmetry]
If $\NF(M_0,M_1,M_2,M_3)$ then
$\NF(M_0,M_2,M_1,M_3)$.
     \item[C7 Disjointness]
If $\NF(M_0,M_1,M_2,M_3)$ then $M_1 \inter M_2 = M_0$.
\end{description}

Axiom C7 is largely a matter of notational convenience.  We will
indicate in Chapter III how the major argument of this paper could be
slightly revised to avoid this axiom.
With it we obtain immediately the following monotonicity property.

\begin{prop}
\label{monotonicity4}
Suppose $NF(M_0,M_1,M_2,M_3)$ and $M_0 \subm N$ which is a {\bf
K}-submodel of $M_1$ and $M_2$.  Then $NF(N,M_1,M_2,M_3)$.
\end{prop}

This result could be easily deduced from
the base extension axiom \ref{baseextension} without using C7.

By taking `formal copies' of $M'_2$ and $M_3$, one derives a variant of
C2 where $M'_1$ can be demanded to be $M_1$.
We refer to Axiom C5 as weak uniqueness because
it simply demands that any two amalgamations of a given vee be
compatible.  Thus, it is making the compatibility class of the diagram
unique, not the amalgamating model.

\begin{lemm}
Suppose
$\langle M_0,M_1,M_2,M_3\rangle$
$\langle M_0,M_1,M'_2,M'_3\rangle$
are full free amalgams.
If
$M'_2$ is
isomorphic to $M_2$ over $M_0$  by a map $g$
then $g$ is an isomorphism between
$M'_2$ and
$M_2$ over $M_1$.
\end{lemm}
\PROOF.  Apply the weak uniqueness
axiom to the map $f$ that is the union of the identity
on $M_1$ and the given $g$ from $M_2$ to $M'_2$.
\endproof

\jtbnumpar{Smoothness}
Does the class {\bf K} admit a `limit' of an ascending chain (or more
generally a directed system) of {\bf K}-stuctures?  There are several
different variants on this question and the answers determine
significant differences in the behavior of {\bf K}.  We discuss the
variants in detail in \cite{BaldwinShelahprimalii}
; we now just mention a couple of
possibilities and some of the consequences.

The strongest requirement is to deal directly with unions of chains.
But even here there are several variations.  One can demand that any
union of a (continuous) increasing chain of {\bf K} stuctures be a
member of {\bf K}. More subtly, one can ask that if each member of the
chain is {\bf K}-submodel of a fixed $M$ then the union is also.

Shelah
has shown that a class {\bf K} satisfying the axioms A0 through A4
enumerated here and stringent requirements for closure under unions can
be presented as the collection of models in a pseudoelementary class
in an infinitary logic
which omit a family of types.  (See Section 1 of
\cite{Shelahnonelemii} and
\cite{Makowskyabstractembed}.)

Beginning in Chapter II we discuss the ways in which unions of chains
can be replaced by demanding the existence of a prime model over the
union.
Again, there are a number of smoothness properties that can be discussed
in this context.  One obvious application is attempts to improve the
L\"owenheim Skolem property to demand that each set can be imbedded in a
model of roughly the same size.
\section{Examples}
 We describe in  this section a number of examples of classes and
notions of amalgamation which satisfy at least some of the axioms that
we are discussing.  Of course, the prototype is the collection of models
of a stable first order theory where a free amalgamation is one that is
independent in the sense of nonforking.

In this section we first discuss some
contrived examples which although they lack any intrinsic interest make
it easy to exhibit some of the pathologies that we are investigating.
Then we place in context some classes which naturally arise in the
attempt to extend classification theory to, e.g., infinitary classes.
\jtbnumpar{Contrived Examples}
\label{contex}
Let {\bf B} be the class of all
structures of the following
sort.  We work in a language with two unary predicates, $U$, $V$ and a
binary relation $E$.  Now $M$ is in {\bf B} if via $E$, each member
$a$ of
$V$ is the name of a subset $X_a =\{m\in U(M):E(m,a)\}$
of $U$ and every subset of $U$ has one and only one
name.  Thus each member $M$ of {\bf B} is determined up to isomorphism
by the cardinality of $U(M)$.
$\NF(M_0,M_1,M_2,M_3)$
holds just if for each $a \in M_1$ ($M_2$), $X_a$
in the sense of $M_3$ is a subset of $M_1$ ($M_2$).
We can illustrate the axioms by defining $\subm_{\bf B}$ in two
different ways.
\begin{enumerate}
\item
Define $M \subm_ {\bf B} N$ if
$U(M) \subseteq U(N)$,
$V(M) \subseteq V(N)$ and each element of $V(M)$ names in $N$ the
same subset of $U(M)$ that it names in $M$.
Under this definition if
$M_0 \subm_{\bf B}
M_1 \subm_{\bf B} M_3$
and
$M_0 \subm_{\bf B}
M_2 \subm_{\bf B} M_3$
then
$\NF(M_0,M_1,M_2,M_3)$.  Moreover if
$M_0 \subm_{\bf B} M_1$
and
$M_0 \subm_{\bf B} M_2$, we can find a common
{\bf K}-extension for them
letting $M_3$ be $M_1 \union M_2$ together with a collection of names
for sets that intersect both $U(M_1)$ and $U(M_2)$.
\item  On the other hand, let $M \subm_{\bf B} N$ just mean that $M$ is
a substructure of $N$.  Now it is still possible to verify Axiom C2.  If
a subset $X$ of $M_0$ is named by elements $a$ of $M_1$ and $b$ of $M_2$
then $a$ and $b$ can be identified by the embeddings into $M_3$.  Using
this strategy to amalgamate
Axiom C7 fails; however the strategy outlined below when we consider
an additional predicate $Q$,
which is needed then
to obtain even C2, will also show that C7 is verified.
\end{enumerate}

Some of the problems with amalgamations become clear if we add a unary
predicate $Q$ and demand that for any subset $W$ of $U(M)$ with power
less than $\kappa$ there exist both $a$ and $b$ in $V(M)$ with $a$
satisfying $Q$ but $b$ not satisfying $Q$ so that $X_a$ and $X_b$
both contain $W$.  (That is, both $Q$ and its complement are dense.)
Axioms C2 and C5 hold under the first definition of $\subm_{\bf B}$.
To make C2 hold under the second interpretation of $\subm$ we must
deal with
a subset $X$ of $M_0$ that has one name in $Q(M_1)$ and another in
$\neg Q(M_2)$.  Now another stategy works.  Add points to the
set attached to one of the names and then fill out the model as freely
as possible.
It is easy to see (cf. Chapter II) that there are
extensions of chains which are incompatible.
\par
${\bf B} _{\kappa}$ is defined in the same way but with the
additional restriction that each $|X_a| < \kappa$.

Another variant arises by replacing the single binary relation $E$ by
a family of binary relations $E_i$ such that for each $i<\kappa$ and
each
$a\in V(M)$, there is a unique $m\in U(M)$ with $R_i(a,m)$.  Thus we
code $\kappa$ sequences rather than sets.

\jtbnumpar{$\aleph_1$-saturated models}
Let $T$ be a strictly
stable first order theory and let $K$ be the class of
$\aleph_1$-saturated models of
$T$.  Take $\NF_{\bf K} $ as nonforking and $\leq_K$ as elementary
submodel.  Then the basic axioms are clearly satisfied but $K$ is not
closed under unions of chains of small cofinality (not fully smooth).
But there are prime models (${\bf F}^{s}_{\aleph_1}$ in the notation of
\cite{Shelahbook} or {\bf SET}$_{\aleph_1}$ in the notation of
\cite{Baldwinbook}) over such chains.
\par
The theory REF$_{\omega}$ of countably many refining equivalence
relations  has {\bf K}-prime models of over chains of cofinality
$\omega$
but they are not minimal.  This leads to $2^{\lambda}$ {\bf K}-models of
power
$\lambda$ when $\lambda^{\omega} = \lambda$.   This argument is treated
briefly in \cite{Shelahmaingap1} (the didip)
and will be reported at more length in the current series of papers.
\par
On the other hand if $T$ is a two dimensional stable theory (cf.
Theorem V.5.8
of
\cite{Shelahbook}) then ${\bf I}(\aleph_{\alpha},{\bf K}) \leq |\alpha +
1|$.
\par
\jtbnumpar{Universal Classes}
See II.2.2 of \cite{Shelahuniversal}.
\jtbnumpar{Finite diagrams stable in power}
See II.2.3 of \cite{Shelahuniversal}.
\jtbnumpar{Infinitary Classes}
See \cite{Shelahnonelemii}.
\jtbnumpar{Banach Spaces}
See \cite{ShelahBS}.
\begin{quest}  Are there are stable universal theories of Banach
Spaces beyond the $L^p$-spaces?
\end{quest}
\chapter{Prime models over diagrams}

In Section 1
we discuss diagrams of models and
the basic properties of prime
models over diagrams.
Section 2 concerns prime models over independent pairs of models.
We also discuss several possibilities for the relation of a
model $M_1$ which is independent from an $M_4$ over $M_0$ with
an $M_2$ with $M_0 \subseteq M_2\subseteq M_4$.
Surprisingly, \ref{transind} deduces a property of the dependence
relation but needs the properties of prime models for the proof.
In Section 3 we discuss prime models over chains.

\section{Diagrams}

We consider two kinds of diagrams.  The first is more abstract because
the partial ordering among the structures is witnessed by {\bf
K}-embeddings.  In the (second)
case of a concrete diagram the partial ordering is
witnessed by actual set theoretic containments.  The existence axioms
for free amalgamations have the general form: Given an abstract diagram,
there is a concrete stable diagram which is isomorphic to it (in the
category of diagrams).

The discussion of abstract diagrams is essential at this stage in
the development of theory.  Once a `monster' model or global
universe of discourse has been posited, one can assume that all
diagrams are concretely realized as subsets of the monster model.
The monster model is easily justified in the first order context
so this subtlety does not seriously arise.  In the more general
case one must worry about it until `smoothness' and thus the
monster model are obtained.  (See \cite{Shelahuniversal} and
\cite{BaldwinShelahprimalii}).

We will assume the existence  of prime models over certain simple
diagrams and certain properties relating prime models and independence.
>From this we will deduce the existence of prime models over more
complicated diagrams.  In a later paper when considering prime model
over still more complicated diagrams (tall trees) we will only be able
to obtain a dichotomy between the existence of prime models and
the existence of many models.

We call a triple $\langle M_0, M_1, M_2\rangle$ with
embeddings of $M_0$ into $M_1$ and $M_2$ a `vee' diagram.
Thus, Axiom C2 asserts the existence of free amalgamations over vee's.

\jtbdef  An abstract {\em {\bf K}-diagram} (indexed by
a partial order
$\langle I,\preceq\rangle$) is a
pair: a sequence of models and
a collection of maps.
The models
$\{M_x:x\in I\}$ and
the maps
$\{f_{xy}:x,y\in I, x \preceq y\}$
must satisfy the following conditions.
\begin{itemize}
   \item $M_x \in {\bf K} $.
   \item If $x \preceq y$ then there is a {\bf K}-embedding $f_{xy}$
from $M_x$ into $M_y$.
    \item  If $x \preceq y \preceq z$ then $f_{xz} = f_{yz} \circ
f_{xy}$.
 \end{itemize}
We will often suppress mention of the maps.
When we need to refer to them we will often describe the family
$\{f_{xy}:x,y\in I\}$
by $f$ and write $f|Y$ for
$\{f_{xy}:x,y\in Y\}$
when $Y$ is a subset of $I$.
\endjtbdef

\jtbdef  $\Mun$ and $\Nun$ are {\em isomorphic {\bf
K}-diagrams} if they are both indexed by the same partial order $I$ and
for each $x\in I$ there is an isomorphism $\alpha_x$ between $M_x$ and
$N_x$ such that the $\alpha_x$ and $f_{xy}$ in $\Mun$ and $\Nun$
commute in the natural way.
We will write $\alpha$ for the family of the $\alpha_x$.

Note that no single function $f$ is a morphism for {\bf K} unless its
domain is in {\bf K}.  In particular, we cannot speak of a single map
whose domain is a vee.

\jtbdef  A concrete
{\bf K}-diagram $\Mun$ inside $M$ (indexed by
$I$) is a collection,
$\Mun =\{M_x:x\in I\}$,
of
members of $K$
such that  each $M_x \subm M$
and if  $x \preceq y$ then $M_x
\subm M_y$.

We extend the notion of a free amalgam to more general diagrams.
\jtbdef  A {\em stable
{\bf K}-diagram} $\Mun$ inside $M$
(indexed by the lower
semilattice
$\langle I,\wedge\rangle$) is a {\underline {concrete}} {\bf K}-diagram
inside $M$ satisfying the additional condition
\begin{itemize}
   \item For any $x,y \in I$, $M_x$ and $M_y$ are freely amalgamated
over $M_{x \wedge y}$ within $M$.
\end{itemize}

\jtbdef
\begin{enumerate}
\item  $N$ is a {\em {\bf K}-extension} of $\Mun$ if there is an
isomorphism $f$
between $\Mun$ and a concrete {\bf K}-diagram $\Nun$ inside $N$.
$f$ is called a {\em {\bf K}-embedding} of $\Mun$ into $N$.

\item  $N$ is a {\em stable
{\bf K}-extension} of $\Mun$ if there is an isomorphism
$f$
between $\Mun$ and a stable {\bf K}-diagram $\Nun$ inside $N$.
$f$ is called a {\em stable {\bf K}-embedding} of $\Mun$ into $N$.
\end{enumerate}

The phrase `$\Mun$ is a stable {\bf K}-diagram' means $\Mun$ is stable
in some $M$ which is suppressed for convenience.
A vee $\langle M_0, M_1, M_2\rangle $ is stable if there
is an $M_3$ completing it to a full free amalgam.
Note that an isomorphism between diagrams need not preserve stability.
It only guarantees isomorphism between the individual $M_x$ and $M'_x$.
In particular an isomorphic copy of a free amalgam need not be free.  So
we cannot say a vee diagram is free without specifying an ambient model
and an embedding into it.
This is true even in the first order case.  In this case we usually
think of the monster model as the ambient model and identity as the
embedding.
\jtbnumpar{Examples}  Of course a (full) free amalgam is a stable {\bf
K}-diagram.
So is any chain of {\bf K}-models which have a common
{\bf K}-extension. Another natural example is an `L'-diagram (see
figure) where for each $i \leq \delta$, $\ind {M_i} {M_0} N$ inside
the ambient model.
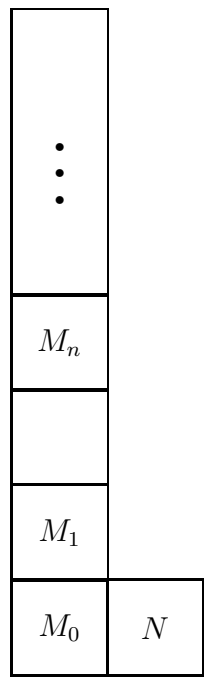
\begin{figure}
\begin{center}
\begin{picture}(76,240)
\put(0,0){\framebox(36,36){$M_0$}}
\put(0,36){\framebox(36,36){$M_1$}}
\put(0,72){\framebox(36,36){}}
\put(0,108){\framebox(36,36){$M_n$}}
\put(0,144){\framebox(36,108){}}
\put(18,180){\circle*{3}}
\put(18,190){\circle*{3}}
\put(18,200){\circle*{3}}
\put(36,0){\framebox(36,36){$N$}}
\end{picture}
\end{center}
\caption{An `L' diagram}
\end{figure}

Compatibility of two extensions of a single model is a straightforward
notion.

\jtbdef  Two members $N_1,N_2$ of $K$ are said to be {\em compatible}
over the
{\bf K}-embeddings
$f_1,f_2$
of $N_0$ into $N_1,N_2$
if there is an $N_3$ in $K$ with both $N_1$ and $N_2$ {\bf K}-embeddible
in $N_3$ by {\bf K}-embeddings $g_1, g_2$ such that
$g_1 \circ f_1
= g_2 \circ f_2$.

In view of the existence axiom for free amalgamations any two members
$M_1$, $M_2$
of $K$ with an $M_0\in {\bf K} $ that is a {\bf
K}-substructure of
each of them are compatible over $M_0$.  When $N_0$ is replaced
by a diagram $\Mun$ involving
infinitely many members
of $K$ the situation is more complicated.

\jtbdef
\label{compatibleembeddings}
Suppose $\Mun$ is a {\bf K}-diagram.
Let $f_1, f_2$ be {\bf K}-embeddings of $\Mun$ into $M_1$ and
$M_2$ respectively.  $M_1$ and $M_2$ are {\em {\bf K}-compatible over
$\Mun$} via $f_1, f_2$
if there  exists an $N\in K$ and {\bf K}-embeddings $g_1, g_2$
of $M_1, M_2$ into $N$ such that $g_1\circ f_1$ and $g_2 \circ f_2$
agree on $\Mun$.

(The $f_i$ are families of maps and the $g_i$ are single maps but the
meaning of the composition should be clear.)


Clearly over each $\Mun$, compatibility defines a reflexive and
symmetric relation.
Axiom C2 shows that the relation is also transitive so for any diagram
$\Mun$ we have an equivalence relation, compatibility over $\Mun$.
\jtbdef
The abstract {\bf K}-diagram
$\Mun$ is {\em stably univalent} if all stable {\bf
K}-extensions of $\Mun$ are compatible.

In this language
Axiom C5 asserts that every  vee is stably
univalent.


There is little hope to find a `prime' model (in the usual categorical
sense) over an arbitrary diagram $\Mun$.  For,
if a diagram isomorphism from $\Mun$ into $N$
collapsed two elements of models $M \in
\Mun$ then there could be no isomorphism from a structure containing
all members of $\Mun$ into $N$.
Thus, we restrict the
following definitions to stable embeddings.
We first define the notion of a `prime' model within a compatibility
class.  Remember that the image of a stable  embedding
is required to be stable.
\jtbdef
\label{compatibilityprime}
\begin{enumerate}
     \item
Let $f$ be a stable  {\bf K}-embedding of the abstract diagram
$\Mun$ in $M$. We say
$M$ is {\em compatibility prime}  over $(\Mun, f)$
(or over $\Mun$ via $f$)
if for every
$M' \in {\bf K} $ and stable embedding $f'$ of $\Mun$ into $M'$
such that $M$ and $M'$ are compatible via $f,f'$
there is an embedding of $g:M\mapsto M'$ with $g \circ f = f'$.
    \item We omit $f$ if it is an identity map.
\end{enumerate}
\par
Strictly speaking, we should refer to a triple $\langle \Mun, f, M
\rangle$.  A little looseness to ease reading seems acceptable here.  We
are more precise  when we introduce the notion of canonically prime
in Section~\ref{primemodelschains}
since a new basic relation is added to the system.

Note that if the diagram $\Mun$ has a unique maximum element then that
element is compatibility prime over $\Mun$.
Note that if $M$ is compatibility prime over $\Mun$ via $f$ then $M$
is isomorphic to an $M'$ which is compatibility prime over $\Mun$ via
the identity.

There are various ways in which a compatibility prime model can fail to
be unique.  There could be more than one compatibility class;  within a
given compatibility class there could be nonisomorphic compatibility
prime models.

\jtbdef
$M$ is {\em absolutely prime} over $(\Mun,f)$ if $M$ can be embedded over
$(\Mun,f)$ into any stable {\bf K}-extension $N$ of $\Mun$.

Clearly if $M$ is compatibility prime over $\Mun$ and $\Mun$ is
univalent then $M$ is absolutely prime.
In view of the weak uniqueness axiom if there is a prime model over a
stable diagram $\Mun$ then $\Mun$ is stably univalent.
We will
see that with strong enough hypotheses on {\bf K} the various notions of
prime coalesce.
We have introduced the notion of `absolutely prime' to emphasize
the distinction with compatibility prime;  `absolutely prime' is the
natural extension of the usual model theoretic notion to `prime over a
diagram'.  When we write `prime' with no adjective we mean `absolutely
prime'.

\section{Prime Models over vees}
\label{primemodelsovervees}
In the next two sections we discuss axioms concerning the existence and
properties of prime models over certain specific diagrams.
We begin by assuming the existence of prime models over vee's.  In the
light of C5 (weak uniqueness) all completions of an amalgam are
compatible
so it would be equivalent to replace absolutely prime by compatibility
prime in the following axiom.  We will often just say `prime' for the
absolutely prime model over a vee.  Again because of weak uniqueness we
don't
have to specify $M_3$ is prime over $M_1 \union M_2$ {\it in $M_4$}.

\jtbnumpar{Prime Models over vees}
\mbox{}
\begin{description}
    \item[D1]  There is an absolutely
prime model over any
free amalgam $\langle M_0, M_1, M_2\rangle$.
\end{description}

For compactness we write $M$ is prime over $M_1 \union M_2$ instead of
$M$ is prime over $\langle M_0, M_1, M_2 \rangle$.

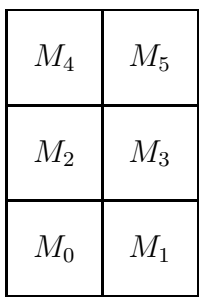
\begin{figure}
\begin{center}
\begin{picture}(76,120)
\put(0,0){\framebox(36,36){$M_0$}}
\put(0,36){\framebox(36,36){$M_2$}}
\put(0,72){\framebox(36,36){$M_4$}}
\put(36,0){\framebox(36,36){$M_1$}}
\put(36,36){\framebox(36,36){$M_3$}}
\put(36,72){\framebox(36,36){$M_5$}}
\end{picture}
\end{center}
\caption{Base Extension Axioms}
\label{basextfig}
\end{figure}

Now we describe properties which relate independence and
prime models.
The next axiom of this group
corresponds to Axiom C4 of \cite{Shelahuniversal}.
Diagram \ref{basextfig}
illustrates each of the principles discussed below.

\jtbnumpar{The base extension axiom}
\label{baseextension}
\mbox{}
\begin{description}
\item[D2]
Suppose $M_1$ and $M_4$ are
freely amalgamated over $M_0$ in $M_5$.  If   $M_0\subm M_2
\subm M_4$ and $M_3$ is prime over $M_1 \union M_2$ in $M_5$
then
$M_3$ and $M_4$ are freely amalgamated over $M_2$ in $M_5$.
\end{description}

To understand the base extension axiom in the first order context, think
of $M_3$ as $M_2[M_1]$ (For this notation see \cite{Baldwinbook}.).
Then the axiom is implied by the fact that for
any $X$ and any model $M$ (in the appropriate category, e.g.
$\omega$-stable and the normal notion of prime model), $X$ dominates
$M[X]$ over $M$.

Axiom~\ref{baseextension}
yields a somewhat surprising consequence.  We can obtain the
following
`transitivity of nonforking' from \ref{baseextension} and the
weak uniqueness
we posited in Axiom C5.  This result is remarkable because we are
establishing a property of the nonforking relation which makes no
reference to prime or generated models.  But, the proof uses properties
of either generated models
(the version in \cite{Shelahuniversal}) or prime models
(here).

\begin{thm}[Transitivity of independence]
\label{transind}
\mbox{}
If $\NF(M_0,M_1,M_2,M_3)$ and $\NF(M_2,M_3,M_4,M_5)$
{then $\NF(M_0,M_1,M_4,M_5)$.}
\end{thm}
\PROOF.
By the existence axiom there are $M''_4$ and $M''_5$ and an isomorphism
$g$ of $M_4$ and $M''_4$ over $M_0$ such that $\ind {M_1} {M_0} {M''_4}$
in $M''_5$.  By monotonicity (C3\k/ iii)) we may assume $M''_5$  is
prime over $M_1 \union M''_4$.  Let $M''_2$ denote $g(M_2)$ so $M_0\subm
M''_2\subm M''_4$.
We want to show the existence of an isomorphism with domain $M''_5$
which fixes $M_1$ and
maps $M''_4$ to $M_4$.

Let $M'_3$ be prime over $M_1 \union M_2$ and contained in
$M_3$.  By the monotonicity axioms (C3\k/ iii) and C3\k/ i)) we have
$\ind
{M_1} {M_0} {M_2}$ in $M'_3$ and $\ind {M'_3} {M_2} {M_4}$ in $M_5$.
Similarly, if $M'_5$ is chosen prime over $M'_3\union M_4$ in $M_5$ then
$M'_5 \subm M_5$ and  $\ind {M'_3} {M_2} {M_4}$ in $M'_5$
using axiom C3\k/ iii).

Since $M_2$ and $M''_2$ are isomorphic over $M_0$ and both are
independent from $M_1$ over $M_0$ there is an isomorphism $f$ taking
$M'_3$ into $M''_5$ which extends $g|M_2 \union 1_{M_1}$.  (This follows
because $M'_3$ is prime over $M_2 \union M_1$ and applying the
weak uniqueness axiom.)  Let $M''_3$ denote $f(M'_3)$.  Then $M''_3$ is
prime over $M''_2 \union M_1$ since $f$ is an isomorphism.  By the base
extension axiom  \ref{baseextension}
$\ind {M''_4} {M''_2} {M''_3}$ in $M''_5$.
Let $M'''_5$ be prime
over $M''_4 \union M''_3$ and let $h$ be a map from $M'''_5$ into some
$N$
which extends $g \union f^{-1}$.  Now $M'_4 = h(M''_4)$ and $M_4$ are
isomorphic
over $M_2$ and both are independent from $M'_3$ over $M_2$.  So by the
weak uniqueness axiom there is an $h'$ defined on $N$ which takes $M'_4$
to $M_4$ and fixes $M'_3$.  Now $h'\circ h$ takes $M''_4$ to $M_4$ and
fixes $M_1$ as required.
\endproof
\par
\begin{thm}[Transitivity of primality]
\label{transprime}
Suppose $M_1$ and $M_4$ are
freely amalgamated over $M_0$ in $M_5$.  If   $M_0\subm M_2
\subm M_4$,
$M_3$ is prime over $M_1 \union M_2$ in $M_5$
and
$M_5$ is prime over $M_3 \union M_4$ in $M_5$
then
$M_5$ is prime over $M_1 \union M_4$ in $M_5$.
\end{thm}
\PROOF.
Let $N$ be a stable {\bf K}-extension of $\langle M_0, M_1, M_4\rangle$.
That is, there exist maps
$f_1$, $f_4$ from
$M_1$, $M_4$
to {\bf K}-submodels
$M'_1$, $M'_4$ of $N$
that agree on $M_0$;  let $M'_0$ denote the common image of the $f_i$
on $M_0$.  Then since the embedding is stable
$\ind
{M'_1} {M'_0} {M'_4}$ inside $N$.
We must find a common extension of the $f_i$ to
$M_5$.  Let $f_2$ denote the restriction of $f_4$ to $M_2$ and $M'_2$
its image.  Now, since $M_3$ is prime over $M_1 \union M_2$ there is a
map $f_3$ with domain $M_3$ which extends $f_1 \union f_2$.  Denote the
image of $f_3$ by $M'_3$.  By the base extension axiom
\ref{baseextension}
$\ind
{M'_3} {M'_2} {M'_4}$ in $N$.  So there is an $f_5$ mapping $M_5$ into
$N$ and extending
$f_3 \union f_4$.
A fortiori, $f_5$ extends
$f_3 \union f_4$ and we finish.
\endproof

There are some further properties of prime models which both arise in
some natural situations and are useful tools.  We describe them now but
they do not play an important role in the theory until we reach some
rather special cases.

\jtbdef  The concrete {\bf K}-extension $M$ of $\Mun$ is {\em minimal}
over the diagram $\Mun$ if there
is no proper {\bf K}-submodel of $M$ that contains  $\Mun$.
\endjtbdef

\jtbnumpar{Some further Axioms}
\mbox{}
\begin{description}
    \item[D3]  The prime model over a free amalgam $\Mun$ is
minimal over $\Mun$.
    \item[D4]  The prime model over a free amalgam $\Mun$ is
unique up to isomorphism over $\Mun$.
\item[D5]
\label{basextb}
Suppose $M_1$ and $M_4$ are freely amalgamated over $M_0$ in $M_5$ and
$M_5$ is prime over $M_1 \union M_4$.
If
$M_0 \subm M_2 \subm M_4$ and $M_3$ is prime over $M_2 \union M_1$ then
$M_5$ is prime over $M_4 \union M_3$.
\end{description}

D5 is a kind of `converse' to Theorem \ref{transprime} that we may
need later.
It will only hold in very restricted cases.  It is true for the notion
of generation in universal classes; it fails for prime models in the
first order case (even $\omega$-stable).
However,
if prime models over vees are minimal (e.g. for a first order theory
without the dimensional order property ) then an
even stronger version of D5 holds easily.  Namely
the monotonicity requirement that if $M_5$ is minimal over $M_1 \union
M_4$ then $M_5$ is minimal over $M_3 \union M_4$.

\section{Prime models over chains}
\label{primemodelschains}

An abstract chain
$\langle M_i, f_{i,j}
:i,j\in I\rangle$
is an abstract diagram
whose index set $I$ is linearly ordered.  Any closed initial segment of
an abstract chain has a natural representation as a concrete chain.  To
see this consider
$\langle M_i, f_{i,j}
:i,j \leq \alpha \rangle$.  Let $M'_i$ denote $f_{i,\alpha}(M_i)$.
Then the $M'_i$ are a concrete chain with inclusion maps
$f'_{i,j} =
f_{j,\alpha} \circ
f_{i ,j} \circ
f^{-1}_{i,\alpha}$.  Thus for any chain indexed by
an ordinal $\gamma$
and any limit ordinal $\delta <\gamma$ it makes sense to speak
of $
\union_{i<\delta}M_i$ as we can concretely realize $\Mun |\delta$ in
$M_{\beta}$ for any $\beta$ with $\delta \leq \beta \leq \gamma$.

We discuss in this section the specification and existence of prime
models over chains.  In general, given an increasing chain of {\bf
K}-models there is no reason to assume that the chain has any common
extension in {\bf K}, let alone one that is prime.
If we assume the existence but not the
weak
uniqueness of compatibility
prime models, there may be two incompatible compatibility
prime
models over the same chain.
(We say {\bf K} is not smooth.)
For this reason, we can not in the most
general case just define `prime' as compatibility prime
and posit that `prime' models exist.
We would need to introduce another axiom asserting that
there is only one compatibility class over any chain.  Justification
of such an axiom is the main point of \cite{BaldwinShelahprimalii}.
However, to reduce the set theoretic hypotheses of that argument
we introduce a new predicate ($\cpr$) with the intuitive meaning,
`$M$ is canonically prime over $\Mun$' and prescribe axioms describing
the behavior of such prime models.

Consider the last example in Paragraph~
\ref{contex}.  We have a unary predicate
$U$ which is dense and codense.  At any limit stage of cofinality
$\kappa$,
we have to decide whether the name of certain $\kappa$ sequences are
in $U$ or not.  Different answers correspond to different compatibility
classes.

We want to demand the existence of a `prime' model over a union of a
chain.  If the chain has length longer than $\omega$ several
possibilities arise for what we should demand of models at limit stages
in the chain.  It is unreasonable to demand the existence of `prime'
models over chains that are not {\bf K}-continuous in the following
sense.

\jtbdef
\begin{enumerate}
\item
The chain $\langle M_i, f_{i,j}
:i,j < \beta\rangle $ is {\em
{\bf K}-continuous}
if for each limit ordinal $\delta < \beta$,
$\cpr(\Mun_{\delta}, M_{\delta},f|\delta)$.
\item
The chain $\langle M_i, f_i
:i < \beta\rangle $ is {\em
continuous}
if for each limit ordinal $\delta < \beta$,
$M_{\delta}= \union_{i<\delta}M_i$.
\end{enumerate}

Axiom Ch1 gives an implicit definition of a canonically prime model;
Axiom Ch2 asserts that such a model exists.  Suppose $\Mun$ is a `long'
chain.  At each limit stage one has a choice of compatibility classes
for a `prime' model over that segment of the chain.  The notion of a
canonically prime model requires that these choices cohere.

\begin{description}
    \item[Ch1]
$\cpr(\Mun,M,f)$ implies
\begin{enumerate}
    \item $\Mun$ is a {\bf K}-continuous chain.
    \item $M$ is compatibility prime over $\Mun$ via $f$.
\end{enumerate}
We often write
$M_\delta$ is
canonically prime over $\Mun_{\delta} = \langle M_i:i < \delta\rangle$
for the formal expression
$\cpr(\Mun_{\delta}, M_{\delta},f|\delta)$.  That is, for brevity we do
not
mention the specific embeddings $f_{i,j}$ unless they play an active
role in the discussion.
    \item[Ch2]  For any {\bf K}-continuous chain $\Mun$ there is
a model $M$ and a family of maps $f$ that satisfy $\cpr(\Mun,M,f)$.
    \item[Ch3]  The canonically
prime model over an increasing chain $\Mun$ is
unique up to isomorphism over $\Mun$.
\end{description}
Again the uniqueness axiom is regarded as a desirable property to prove
and is not assumed in the general development.
\jtbdef The chain $\langle M_i, f_{i,j} :i < j < \beta\rangle $
is called {\em essentially
{\bf K}-continuous}
if for each limit ordinal $\delta < \beta$,
there is a model $M'_{\delta}$ which can be interpolated
between the predecessors of $M_{\delta}$ and $M_{\delta}$ and
is canonically prime over $\Mun_{\delta}$.

More formally, we add to the index set a new $\delta'$ for
each limit ordinal $\delta$.
There exists a system of embeddings $g$
such that if $\alpha$ and $\beta$ are successor ordinals
$g_{\alpha, \beta} =
f_{\alpha, \beta}$,
and for each limit ordinal $\delta$ and each $\alpha < \delta$
$g_{\alpha,\delta}$
factors as
$g_{\alpha,\delta'}\ocirc
g_{\delta',\delta}$.  Finally
$\cpr(\Mun_{\delta}, M_{\delta'},g|\delta')$.

The next axiom provides a local character for dependence of the prime
model over a chain over another model.  This is our only basic axiom
connecting the independence of an infinite diagram with that of its
constituents.
\jtbnumpar
{Axiom L1: Local Dependence for L Diagrams}
If {$\langle M_i:i \leq \delta\rangle$} is a {\bf K}-continuous
increasing sequence inside $M'$
and for each $i < \delta$, $\ind {M_i} {M_0} {N}$
in  $M'$  then $\ind {M_{\delta}} {M_0} {N}$ in $M'$.
\endjtbnumpar
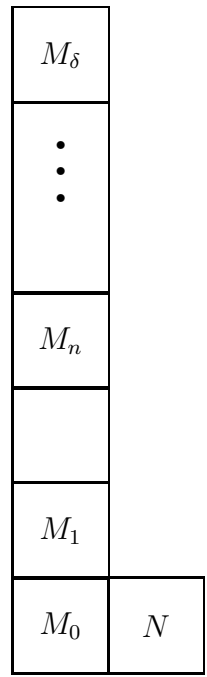
\begin{figure}
\begin{center}
\begin{picture}(76,240)
\put(0,0){\framebox(36,36){$M_0$}}
\put(0,36){\framebox(36,36){$M_1$}}
\put(0,72){\framebox(36,36){}}
\put(0,108){\framebox(36,36){$M_n$}}
\put(0,144){\framebox(36,72){}}
\put(0,216){\framebox(36,36){$M_{\delta}$}}
\put(18,180){\circle*{3}}
\put(18,190){\circle*{3}}
\put(18,200){\circle*{3}}
\put(36,0){\framebox(36,36){$N$}}
\end{picture}
\end{center}
\caption
{Local Dependence for L's}
\end{figure}

\chapter{Prime models over small loose trees}

In this chapter we extend the existence of prime models over simple
diagrams to obtain prime models over more complicated diagrams.  Ideally
we would show the existence of prime models over an arbitrary
independent tree.  This project runs into difficulties when
considering trees of height greater than $\omega$; as a substitute we
show how to obtain prime models over `loose trees' of models indexed
by subsets of $\lambda^{<\omega}$.  In a later paper we expect to
reduce the
discussion of prime models over trees with large height to treatment of
loose trees of height $\omega$.


While the original motivation for loose trees was the reduction of
problems about tall trees to problems about loose trees of countable
height, it turns out that loose trees are closed under several useful
operations such as quotient
and trees are not.  Thus, by passing through loose trees we
are able to obtain results about trees that, at least a priori, are
otherwise unavailable.

Consider a concrete stable diagram.  An initial strategy for building a
prime model over $\Mun$ is to
enumerate $\Mun$ say as $M_i$ for $i < \alpha$
and choose a family of models $N_i$ for $i<\alpha$ so that for each $i$,
$M_{i+1}$ is independent from $N_i$ over the predecessors of $M_{i+1}$
in the tree and then to take $N_{i+1}$ prime over $M_{i+1}$ and $N_i$.
Take canonically prime models over the earlier $N_i$ at limits.  The
resulting model clearly depends on the order of enumeration.  Can one
still prove that this model is compatibility prime over the diagram?  We
show that the answer is yes
if the diagram is a `short' tree (and generalize to
allow loose trees).  Finite trees are considered in 
Section~\ref{freeloosetrees}; trees of countable height in Section
\ref{Locallyfreeloosetrees}.  In the case of finite trees
we show fairly directly that if a loose tree is free under one
enumeration then it is free under any enumeration.  In the second case
we pass to the ostensibly more general notion of locally free loose
tree, show the existence of prime models over such a tree, and deduce
from that the fact that a locally free loose tree is free under any
enumeration.

{\bf Assumption:  An Adequate Class.}
\label{adequate}
We assume in this chapter axiom groups A and C,
axiom
D1 and D2 from group D,
Ch1, Ch2 from group Ch,
and (beginning with
\ref{primoverdiag}) L1.
We call a class with these properties an {\em adequate class}.
\section{Free Loose Trees}
\label{freeloosetrees}
In this section we show that if $\Mun$ is a {\em finite} free loose tree
(definitions follow) of models from an adequate class
then there is an absolutely prime model over
$\Mun$.  More precisely, we say
that $N$ is {\em explicitly prime}
over a finite loose tree $\Mun$ if it is the last in a sequence of prime
models over vee's satisfying certain conditions.  We show that
$N$ is explicitly prime
over a subdiagram
of $\Mun$ then
the sequence witnessing this can be extended to one witnessing
the existence of
an explicitly prime model over $\Mun$.
This ostensibly technical result is essential for the discussion of
prime models over {\em infinite} loose trees in the next section.

\jtbnot A tree
{\bf T} is a partially ordered set which is isomorphic to a
subset of
$\lambda^{<\omega}$  which is closed under initial segment.  We will
often deal directly with this representation.
A tree is
partially ordered by containment.  If $s,t \in {\bf T} $ then $s \wedge
t$ denotes the largest common initial segment of $s$ and $t$ and for any
$t$ other than the root, denoted $\emptynode$, $t^-$ denotes the
predecessor of $t$.

As in the study of stable diagrams we work only with embeddings that
are inclusions, that is, with concrete diagrams.
\jtbdef  A {\em loose tree of models}
$\Mun =\{M_t:t\in {\bf T} \}$ inside $M$
indexed by a tree
{\bf T} is a collection of {\bf K}-models such that if
 $t^- = s$ then
 $M_t \inter M_s \subm M_t$
and
$M_t \inter M_s \subm M_s$.

Note that any stable diagram indexed by a tree
is a loose tree.

We have {\em not} introduced a loose tree as a kind of partial
order.  All our index sets are trees in the normal sense.  A loose tree
{\em of models} is loose because of the inclusions amongst the
models determined by the indexing.

In defining an independent loose tree (below) we will speak of
independence over $M_t \inter M_{t^-}$.  This use of intersection
depends on
Axiom C7
(which implies that $M_s \inter M_t \in {\bf K} $).  However, at a cost
in complexity of notation the results of this chapter could be obtained
through
modifying the definition of a loose tree and a free loose tree
by replacing $M_t \inter M_s$ by a substructure $M_{t,s} \in {\bf K} $.
This would allow us to extend the definition of loose trees to abstract
diagrams.  With such an extension the analogy
would be closer
between the notion
of compatibility prime (over a stable embedding of an abstract
diagram \ref{compatibilityprime})
and the definition
\ref{defprimloosetree}
of a compatibility prime model over a loose tree.

\jtbdef  An {\em isomorphism of loose trees} $f:\Mun \mapsto \Mun'$ is
a family of isomorphisms
$f_t$ taking $M_t \in \Mun$ to $M'_t \in \Mun'$ such that for $t
\leq s$, $f_t|(M_t\inter M_s) \subseteq f_s$.  We say the loose tree
$\Mun$ can be {\bf K}-embedded in $N$ if there is such an isomorphism
between $\Mun$ and a loose tree $\Mun'$ inside $N$.
\label{embltree}
\endjtbdef
\jtbdef
A wellordering $\tbar
= \{t_i:i < \beta
\}$, where $|\beta| =
|{\bf T} |$,
such that if $t_i$ precedes $t_j$ in {\bf T}
then $i \leq j$ is called
an {\em enumeration} of the tree {\bf T}.

An enumeration of {\bf T} induces an enumeration of any $\Mun$ indexed
by {\bf T}.
\jtbnumpar{Example}
\label{exampleloosetree}
The following observation does not figure in our argument but
illustrates one of the subtle differences between a loose tree and a
tree of models.  The partial order of a tree of models is determined
by containments among the models.  The partial order of a loose tree
is artificially imposed (and thus cannot be ignored).
Consider a collection of models $M_0 \subseteq M_1$,
$M_0 \subseteq M_2 \subseteq M_3 \subseteq M_4$ with $\ind {M_1} {M_0}
{M_4}$ and $M_1 \inter M_4 =M_0$.  As a loose tree of models they
can be indexed by the integers less than 5 with the following tree
order:  0 is the root, 1 and 2 are incomparable successors of 0, 3 and 4
are {\em incomparable } successors of 2.  Clearly
$\langle
0,1,2,3,4\rangle$
is an enumeration of the loose tree;  so is
$\langle
0,1,2,4,3\rangle$

Again, we could index these models by a tree with $0$ as an initial
element, $2$, $3$, $4$ as incomparable successors and $1$ above $2$.
Now one can enumerate the tree as $\langle 0,2,3,4,1\rangle$.

With this housekeeping out of the way we can introduce a more important
concept.

\jtbdef
\label{freeloose}
  The loose tree $\Mun$ indexed by {\bf T} is {\em free} or {\em
independent} in $M$
with respect to an enumeration $\tbar$ of {\bf T} if there exists
an ordinal $\beta$ and a
sequence of {\bf K}-models $\langle N_i: i < \beta \rangle$ such
that \begin{enumerate}
    \item All the $M_t, N_i \subm M$.
    \item  $M_{t_0} = N_0$,
    \item
Fix $k < i$ such that
$t_i^- = t_k$
\begin{itemize}
          \item If $1\leq i < \omega$
then $\ind
{M_{t_{i}}} {M_{t_{i}}\inter M_{t_k}} {N_{i-1}}$ in $N_{i}$.
          \item If $i \geq \omega$
then $\ind
{M_{t_{i}}} {M_{t_{i}}\inter M_{t_k}} {N_{i}}$ in $N_{i+1}$.
\end{itemize}
    \item  $\langle N_i: i < \beta\rangle$ is {\bf K}-continuous.
\end{enumerate}
Thus $|\beta| = |{\bf T} |$ but when {\bf T} is infinite the ordinality
of $\beta$ will be $\ord (\tbar) + 1$ which
may be greater than $|{\bf T}  |$.
If $\Mun$ is actually a tree the $M_{t_i} \inter M_{t_k}$ in
condition iii) becomes $M_{t_k}$.
We call $\Nun$ a witnessing sequence for the freeness.

\jtbnumpar{Remark}  The two conditions in iii) of \ref{freeloose} could
be combined if we indexed the $N_i$ by $1 \leq i \leq \beta$.  We didn't
make the change to avoid introducing errors in the proofs of the later
theorems using finite trees indexed in accordance with the official
definition.

A trivial induction shows the following refinement of the definition.
We use this observation without comment below.

\begin{lemm} If $\Mun$ is a loose tree in $M$ indexed by the finite
tree {\bf T}, the witnessing sequence $\langle N_i:i < \beta\rangle$
can be chosen to satisfy
$N_{i}$ is prime over
${M_{t_i}}\union N_{i-1}$ inside $N$.
\end{lemm}

With this lemma in mind,
suppose  the loose tree $\Mun$ indexed by {\bf T} is {\em free} or {\em
independent} in $M$
with respect to the enumeration $\tbar$ of {\bf T} and
there is an initial sequence $\langle
M_{t_i}
:i <j\rangle$ of the
enumeration with
$M_{t_i}    \subseteq M_{t_{i+1}}$.  Then for $i < j$, $N_i$ can
be chosen as
$M_{t_i}$.  This is a likely possibility for a tree, an unlikely
one for a properly loose tree of models.
\jtbnumpar{Some generalizations}  The following variant may turn out to
be necessary.  Define a loose tree to be {\em almost free in $N$} with
respect to an enumeration $\tbar$ if it is free in some extension $N'$
of $N$.  Then a compatibility class of $\Mun$ is almost free with
respect
to $\tbar$ if some $\Mun$ is free with respect to $\tbar$ in some $N$ in
the compatibility class.  We will finally show in this section that
almost free implies free by showing the explictly prime model over
$\Mun$ inside $N'$ can be chosen inside $N$.  This paragraph is
analogous to the monotonicity conditions on freeness over a vee.
\par
The following lemma is proved by a straightforward induction on the
length of the enumeration using the existence of prime models over
independent pairs and canonically prime models over chains.
\begin{lemm}
\label{existfreeenum}
Let $\Mun$ be a loose tree in $M$ enumerated by $\tbar$.  There
is an $N$ and an
isomorphic
copy of $\Mun$ inside $N$ that is free with respect to the
enumeration $\tbar$ in $N$.
\end{lemm}

Note that the isomorphism in Lemma~\ref{existfreeenum} is an isomorphism
of loose trees.  Example~\ref{exampleloosetree} is not free in either
enumeration.
We want to replace `free with respect to an enumeration'
in Lemma~\ref{existfreeenum}
by `free'; that
is to show that if $\Mun$ is free with respect to one enumeration it is
free with respect to any enumeration.  Then we will show that
compatibility
prime
models exist over free trees.  We first handle the case of finite trees.
For this we need some combinatorial tools to reduce the main result to
a more
manageable case.
\jtbdef  Let $\sbar$ and $\tbar$ be enumerations of {\bf T}.  Then,
$\sbar$ and $\tbar$ are {\it close neighbors} if for some $i$,
$s_i =
t_{i+1}$,
$t_i =
s_{i+1}$, and they agree on all other arguments.  $\tbar$ and
$\sbar$ are {\it neighbors} if there is a sequence $\tbar =\tbar^0,
\tbar^1, \ldots \tbar^k = \sbar$ such that for each $j<k$, $\tbar^j$
and $\tbar^{j+1}$ are close neighbors.
\endjtbdef
The next lemma reduces many problems about the relations between two
enumerations to the relation between neighbors and thus by easy induction
to relations between
close neighbors.
\begin{lemm}  If $\sbar$ and $\tbar$ are enumerations of {\bf T} then
they are neighbors.
\label{neighbors}
\end{lemm}
\proof.  Fix $\sbar$ and $\tbar$.  Choose $k$ maximal such that for some
neighbor $\tbar'$ of $\tbar$, $\sbar|k = \tbar'|k$.  Choose the least
$l$, necessarily greater than $k$, such that for some such $\tbar'$,
$s_{k+1} = t'_l$.
Note
that
$s_{k+1} = t'_l$ is incomparable with $t'_{l-1}$, $l-1 \geq k$.
Now if we manufacture $\tbar''$ from $\tbar'$ by
switching $t'_l$ and $t'_{l-1}$ we contradict the minimality of
$l$.
\jtbnumpar{Remark}  Our notion of enumeration corresponds to the concept
of `linear extension' in the theory of partial orderings.
We are not aware if this result has been proved in that context although
it seems likely.
\endjtbnumpar
 \begin{lemm}
\label{finind}
If {\bf T} is a finite tree and $\Mun$ is a loose
tree indexed by {\bf T} which is free
inside $N$
with respect to
the enumeration $\tbar$ then $\Mun$ is free inside $N$
with respect to any
enumeration $\sbar$ of $\Mun$.
\end{lemm}
\PROOF.
By Lemma~\ref{neighbors}
it suffices to show that if $\Mun$ is independent
with respect to $\tbar$ then $\Mun$ is independent with respect to
the enumeration
$\tbar\circ (i,i+1)$.  Let $i = k+1$.  Since $\tbar$ and
$s =\tbar\circ(i,i+1)$
are enumerations,
${t_i}$ and ${t_{i+1}}$ are incomparable.
We will construct a witnessing sequence $N'_i$ for the independence of
$\Mun$ with respect to $\sbar$ from the given sequence $N_i$ witnessing
the independence of $\Mun$ with respect to $\tbar$.  For $l \leq k$, let
$N'_l = N_l$.  To simplify notation, let
$M^a$ denote $M_{t_i}$,
$M_a$ denote $M_{t^{-}_i}$,
$M^b$ denote $M_{t_{i+1}}$, and
$M_b$ denote $M_{t^{-}_{i+1}}$.
We have
$\ind {M^a} {M^a \inter M_a} {N_{i-1}}$ in $N_{i}$
and
$\ind {M^b} {M^b \inter M_b} {N_{i}}$ in $N_{i+1}$.
Since $M_b = M_{t_l}$ for some $l < i$, $M_b \subseteq N_{i-1}$.  By
monotonicity, $\ind {M^b} {M_b\inter M^b} {N_{i-1}}$ in $N_{i}$.
Now choose $N'_{i}\subm N_{i+1}$ prime over $M^b \union N_{i-1}$.  By
the base extension axiom $\ind {N_{i}} {N_{i-1}} {N'_{i}}$ in $N_{i+1}$.
Since we also have
$\ind {M^a} {M_a\inter M^a} {N_{i-1}}$ in $N_{i}$
and $M^a \subseteq N_{i}$
transitivity of independence, Theorem \ref{transind}, yields
$\ind {M^a} {M_a\inter M^a} {N'_{i}}$ in $N_{i+1}$.  Since $N'_i \subm
N_{i+1}$ and if $l > i$,
$\ind {M_{t_l}} {M_{t_l^-}} {N_i}$ implies
$\ind {M_{t_l}} {M_{t_l^-}} {N'_i}$
 we can let $N'_l =
N_l$ for $l > i+1$.
\endproof

\jtbdef
Henceforth we will say a finite loose tree of models is {\em free}
if it is
free under some enumeration.
\endjtbdef

We want to find an analogous result for
certain
loose trees indexed by subsets of $\lambda^{<\omega}$.  This requires
several further concepts.
We begin by extending our notion of prime over a stable diagram to free
loose trees.  Compare the following definition to Definitions
\ref{compatibleembeddings} and \ref{compatibilityprime}.
\jtbdef
\label{defprimloosetree}
The model $M$ is {\it compatibility prime over the
free
loose tree
$\Mun$} with respect to the enumeration $\tbar$
 via the
embedding $f$ of $\Mun$ into $M$ if the image of
$\Mun$ under $f$ is free inside $M$ with respect to $\tbar$
and for every $N \in K$
and every
$g$ embedding $\Mun$ into $N$
with the image of $g$ free inside $N$
such that $M$ and $N$ are compatible over $\Mun$ via $f$ and $g$,
there is a $\hat g$ mapping $M$ into $N$
such that $g = \hat g \circ f$.
\endjtbdef
\medskip
As in the definition of prime over a stable diagram, we need some
`freeness' condition on a loose tree before it makes sense to speak of a
`prime' model over it.
If we strengthen the definition by allowing $N$ to be arbitrary
rather than
requiring that $N$ be compatible over $\Mun$ with $M$
we call
$M$ absolutely prime over $\Mun$.

We need a new `categorical' definition for prime over a loose tree
because the hidden requirements on composition of maps are much looser
when a family is indexed as a loose tree than in our notion of prime
model over a stable diagram.

When we restrict to finite loose trees of models we can define
absolutely prime models.
\begin{lemm}
\label{finprime}
If {\bf T} is a finite tree  and $\Mun$ is a loose tree of models
indexed by {\bf T} which is free inside $N$ with respect to an
enumeration $\tbar = \langle t_i:i < k\rangle$
then there is an absolutely prime model over the loose tree
$\Mun$.
\end{lemm}
\PROOF.  An easy induction on
$|{\bf T}|$ shows
$N_{k-1}$ is the required model.
\endproof
We christen the result of the last construction.
\jtbdef
\begin{enumerate}
\item
If $\langle N_i:i<k\rangle$ witnesses that the finite loose
tree
$\Mun$ is free then we say $N_{k-1}$ is {\em explicitly prime} over
$\Mun$ (with respect to the embedding map and a specific enumeration).
\item
If
 $\langle N_i:i<\beta\rangle$
witnesses that the loose
tree
$\Mun$ is free and $\beta$ is a limit ordinal
then we say $N_{\beta}$ is {\em explicitly prime} over
$\Mun$ (with respect to the embedding map and a specific enumeration) if
$N_{\beta}$ is canonically prime over
 $\langle N_i:i<\beta\rangle$.
\end{enumerate}

\endjtbdef
Now we invoke the prime models to strengthen the sense in which
`freeness' is independent of the enumeration of $\Mun$.
In order to state the result we need several further notations.

\jtbdef
\label{treenot}
\begin{enumerate}
\item
An {\em ideal} of a tree {\bf T} is a subset ${\bf T}_1$ of
{\bf T} that is closed under initial segment.
\item
If $I$ is an ideal of {\bf T} then ${\bf T}_I$ denotes the (quotient)
tree whose
elements are $({\bf T}-I)\union \{\emptynode\}$
with the same meets as {\bf T} if possible but
, if
$x\wedge y$
in the sense of {\bf T}
is in $I$,
then
$x\wedge y = \emptynode$
in the sense of
${\bf T} _I$.
\end{enumerate}
\endjtbdef

Any diagram can naturally be restricted to
a subset $X$ of its index set {\bf T}
by forgetting the models attached elements  ${\bf T}-X$.
We describe several conditions
when the restriction $\Mun|I$ of a (loose) tree $\Mun$ is a (loose)
tree.
For loose trees (but not trees) there is a
natural complement or quotient structure
$\Mun_I$ to $M|I$ which we describe in iii) of the next
definition.
\jtbnot
\begin{enumerate}
\item
If $X \subseteq {\bf T}$ then for any set of models $\Mun$ indexed by
{\bf T}
there is a natural notion of the restriction
$\Mun|X$ of $\Mun$ to $I$:
$\Mun|X = \{M_t:t \in X\}$.
\item
$M_I$ denotes some compatibility prime model over $\Mun|I$.
\item
Suppose
$\Mun=\langle M_t:t \in {\bf T} \rangle$
is a free
loose tree inside $M$.
If $I$ is an ideal of {\bf T} then $\Mun_I$
denotes
$\langle M_t:t \in {\bf T}_I \rangle$.
\end{enumerate}

The term
$M_I$ is an abuse of notation since we are choosing one of a number of
possible compatibility prime models.
If $I$ is an ideal of {\bf T} and $\Mun$ is a free
loose tree of models
indexed by {\bf T} it is easy to see that the quotient tree is also free.
Formally,

\begin{lemm}
If
$\Mun$
is a free
loose tree inside $M$ and
$I$ is an ideal of {\bf T} then $\Mun_I$
is a free
loose tree inside $M$ indexed by ${\bf T}_I$.
\end{lemm}

Now we consider substitutions for a model in a loose tree.
\jtbnot
\label{notformainlemma}
Suppose
$\Mun=\langle M_t:t \in {\bf T} \rangle$
is an indexed family of models.
Then $\Mun(N/s)$ denotes the
indexed family of models.
obtained
by replacing $M_s$ by $N$.
$\Mun(N_1/s_1, N_2/s_2)$ is the natural extension of this notation to
allow two substitutions.

We are interested in a number of substitutions in indexed trees.
Most
will be approximations to the result of replacing the root in
the tree $\Mun_I$ by the prime model over $I$:
$\Mun_I(M_I/\emptynode)$.
The next lemma describes two slightly less
simple ways of deriving new loose trees from a given one.  The fact that
loose trees are obtained by the constructions is an immediate
verification;  the more important fact that the second construction
preserves freeness is proved in Lemma~\ref{omissionlemma}.

\begin{lemm}
Suppose
$\Mun=\langle M_t:t \in {\bf T} \rangle$
is a
loose tree inside $M$.
\begin{enumerate}
\item
Fix $s \in {\bf T} $ and  a model $N\subm M$
such that
for each $t$ with $t^- = s$ or $s^- =t$,
$N \inter M_t \subm M_t$ and
$N \inter M_t \subm N$.
Then
$\Mun' = \Mun (N/s)$ is a loose tree inside $M$.
\item
Fix
$s,r\in {\bf T} $ with $r^- =s$ and $M_s
\subm M_r$.  Let ${\bf T}'$ denote ${\bf T} - \{r\}$.
Then
$\Mun' = \Mun |{\bf T} '(M_r/s)$ is a loose tree inside $M$.
\end{enumerate}
\end{lemm}

We have established the notation to state the following lemma.

\begin{lemm}[The Omission Lemma]
\label{omissionlemma}
Suppose $\Mun = \langle M_t: t\in {\bf
T} \rangle$ is finite and
free loose tree inside $M$,
$s,r\in {\bf T} $ with $r^- =s$ and $M_s
\subm M_r$.  Let ${\bf T}'$ denote ${\bf T} - \{r\}$.  Then
$\Mun' = \Mun |{\bf T} '(M_r/s)$ is free.
\end{lemm}
\proof.  Let $\tbar$ enumerate ${\bf T} $ with $r = t_{j+1}$ and $s =
t_j$.  Suppose $\Nun$ witnesses the freedom of $\Mun$ and $|{\bf T}
|=k$.  Now define
\[ t'_i = \left \{ \begin{array}{ll}
                    t_i & \mbox{if $i < j$} \\
                    t_{i+1} & \mbox{if $j\leq i<k-1$}
                   \end{array}
          \right. \]
\[ N'_i = \left \{ \begin{array}{ll}
                    N_i & \mbox{if $i < j-1$} \\
                    N_{i+1} & \mbox{if $j-1\leq i<k-1$}
                   \end{array}
          \right. \]
It is straightforward that $\Nun'$ witnesses the freeness of $\Mun'$
with respect to the enumeration $\tbar'$.
\endproof
The following theorem asserts that the freeness of a finite loose tree
is independent not only of the enumeration but for a given enumeration
of the choice of the prime models $N_i$ witnessing the freeness.
This theorem allows us to
decompose $\Mun$ into $\Mun|I$ and $\Mun_I$.  An advantage of loose
trees in performing this construction is that they allow the mixing of
the models in the tree and the witnessing sequence.
\begin{thm}
\label{rut}
Suppose $\Mun$ is a free loose tree of models inside $M$
indexed by the
finite tree {\bf T} and $I$ is an ideal of {\bf T}.  Suppose also
that
$N$ is explicitly prime over $\Mun|I$ for an enumeration $\tbar$ of $I$.
Then $\Mun_I(N/\emptynode)$ is free inside $M$.
\end{thm}

\proof.  Suppose $|I| =j$.
Extend the given enumeration of $I$ to an enumeration of {\bf
T} and also denote the extension by $\tbar$.  Let $\Nun$ witness
that $\Mun|I$ is free.    Noting that any initial segment of an
enumeration is an ideal, let ${\bf T} _i$ be the  ideal
composed of the first $i$ elements in the enumeration.  Let $\Mun^i =
\Mun_{{\bf T} _i} (N_i/\emptynode)$.  We finish by showing by induction
on $i < |I|$ that $\Mun^i$ is free inside $M$.

If $i =0$, $\Mun^i = \Mun$ and the result is clear.  Suppose we
know $\Mun^i$ is free and consider $\Mun^{i+1}$.  In ${\bf T}_{{\bf T}
_i}$, $t_{i+1}$ is an immediate successor of $\emptynode$.
Lemma
\ref{techenum} below
shows $\Mun^i(N_{i+1}/t_{i+1})$ is free inside $M$.
This implies $\Mun^{i+1}$ is free by the omission lemma and the
induction hypothesis.
Thus, it remains only to prove
Lemma~\ref{techenum}.
\endproof

Note that any finite tree can be enumerated so that
$M_{t^-_{i+1}}\subm M_{t_i}$
for any $i$.

\begin{lemm}
\label{techenum}
Suppose $\Mun$ is a free loose tree of models inside $M$
indexed by the
finite tree {\bf T}.
Let $N \subm M$ and suppose $N$ is prime over $M_v \union M_{v-}$ for
some
$v\in {\bf T}$.  Then $\Mun' = \Mun(N/v)$ is free inside $M$.
\end{lemm}
\proof.  Fix an enumeration $\tbar$ of $\Mun$ with  $t_j = v$.  Suppose
that $\Nun = \langle N_i:i < |{\bf T} |\rangle$ witnesses the freedom of
$\Mun$.  Before constructing the sequence $N'_i$ which will witness the
freedom of the new tree we need an auxiliary sequence.  We define for $j
\leq i < |{\bf T} |$ an increasing chain of models $N^*_i$ beginning
with $N^*_j =M_v = M_{t_j}$ such that
\begin{enumerate}
\item each $N^*_i \subm N_i$
\item
$\ind {M_{t_{i+1}}} {M_{t^-_{i+1}}\inter M_{t_{i+1}}}
{N^*_i}$ inside $N^*_{i+1}$,
\item
$\ind {N^*_{i+1}} {N^*_i} {N_i}$ inside
$N_{i+1}$,
\item
$N^*_{i+1}$ is prime over
$ {M_{t_{i+1}}} \union {N^*_i} $.
\end{enumerate}

{\bf Base Step:}
First note setting $N^*_{j} = M_{t_j}$ satisfies conditions
ii) and iii).
Choosing $N^*_{j+1}$ to satisfy iv), we
must verify iii).  This follows by the base extension axiom applied
 as in Figure
\ref{basestep}.
Recall $M_{t_j} = M_v = N^*_j$.
\begin{figure}
\begin{center}
\begin{picture}(76,120)
\put(0,0){\framebox(72,36){$M_{t_{j+1}} \inter M_v$}}
\put(0,36){\framebox(72,36){$M_v$}}
\put(0,72){\framebox(72,36){$N_j$}}
\put(72,0){\framebox(36,36){$M_{t_{j+1}}$}}
\put(72,36){\framebox(36,36){$N^*_{j+1}$}}
\put(72,72){\framebox(36,36){$N_{j+1}$}}
\end{picture}
\end{center}
\caption{The base step}
\label{basestep}
\end{figure}

{\bf Induction Step:}  Suppose that for some $k \geq j$ we have chosen
$N^*_{k}$ to satisfy ii) and iii) (with $i = k-1$).  Choose
$N^*_{k+1}$
inside $N_{k+1}$
prime over
$M_{t_{k+1}} \union N^*_{k}$ to satisfy i) and iv).  Now
$N^*_{k+1}$ satisfies ii) by monotonicity and iii):
$$\ind {N^*_{k+1}} {N^*_{k}} {N_{k}} \hbox{\rm\ in }
 N_{k+1}$$
follows
from the base extension axiom as in the induction step diagram,
Figure~\ref{inductionstep}.
\begin{figure}
\begin{center}
\begin{picture}(76,120)
\put(0,0){\framebox(72,36){$M_{t_{k+1}}\inter M_{t^-_{k+1}}$}}
\put(0,36){\framebox(72,36){$N^*_{k}$}}
\put(0,72){\framebox(72,36){$N_{k}$}}
\put(72,0){\framebox(36,36){$M_{t_{k+1}}$}}
\put(72,36){\framebox(36,36){$N^*_{k+1}$}}
\put(72,72){\framebox(36,36){$N_{k+1}$}}
\end{picture}
\end{center}
\caption{The induction step}
\label{inductionstep}
\end{figure}

This completes the construction of the $N^*_i$.

Let $u$ denote $v^-$;
note that since
$N^*_j
= M_v$ we have
$\ind {N_{j-1}} {M_v \inter M_u}
{N^*_j}$.
With this as the base an easy induction on $l$ for $j \leq l
\leq k$ shows
$\ind {N_{j-1}} {M_v \inter M_u}
{N^*_l}$.  We will rely on the case $l =k$.

We have from the free enumeration that
$\ind {M_v} {M_u \inter M_v} {N_{j-1}}$.  The base extension axiom then
gives,
since $N$ is prime over $M_v \union M_{v^-}$,
$\ind {N} {M_v } {N_{j-1}}$.
Choose $M'_v\subm N_k$ prime over
$N\union N_{j-1}$.
By transitivity of primeness, Lemma~\ref{transprime},
$M'_v$ is prime over
$M_v\union N_{j-1}$.

>From the base extension axiom and the conclusions of the last two
paragraphs, we deduce
$\ind {M'_v} {M_v} {N^*_k}$. This allows us to perform the $i = j$
step of the following construction.  Choose $N'_i$ so that
\[ N'_i = \left \{ \begin{array}{ll}
                    N_i & \mbox{if $i < j$} \\
                    M'_v & \mbox{if $i = j$}
                   \end{array}
          \right. \]
and for $j < i$ so that $N'_i\subm M$ is prime over $N^*_i
\union
N'_{i-1}$ and $\ind {N'_i} {N^*_i} {N^*_k}$.  The choice of the
$N'_i$ for $i > j$ is a straightforward induction.

Now to complete the proof we must observe that the $N'_i$ witness the
freeness of $\Mun'$.  That is we must show that
$N'_{i+1}$ is prime
over $M_{t_{i+1}}\union N'_i$.
This follows from transitivity of primeness, Lemma \ref{transprime},
since for each $i$ we have
$N^*_{i+1}$ is prime over $M_{t_{i+1}}\union N^*_i$
and
$N'_{i+1}$ is prime over $N^*_{i+1} \union N'_i$.
\endproof
We restate Theorem~
\ref{rut} in a more applicable form.
\begin{cor}
\label{rutcor}
Suppose $\Mun$ is a free loose tree of models inside $M$
indexed by the
finite tree {\bf T} and $I$ is an ideal of {\bf T}.  If $\tbar$ is
an enumeration of $I$ and $\Nun$ witnesses that $\Mun|I$ is free inside
$M$ then $\Nun$ can be extended to a sequence witnessing that $\Mun$ is
free inside $M$.
\end{cor}

We have shown that if {\bf K} is an adequate class and $\Mun$ is a
finite (loose) tree $\Mun$ which is free inside $N$ under some
enumeration then $\Mun$ is free under any enumeration and there is an
absolutely prime model over $\Mun$.  With more difficulty we showed that
if $I$ is an ideal of $\Mun$ and $N_0$ is explicitly prime over $I$ then
there is an sequence defining an explicitly prime model over $\Mun$
which includes $N_0$.
\section
{Locally Free Loose Trees}
\label{Locallyfreeloosetrees}
We now extend our analysis to infinite trees.  Recall that in the
absence of the monster model we can only speak of a diagram being free
when we have an embedding into an ambient model in mind.  And only then
can we discuss the possibility of any
sort of `prime' model over the
diagram.  This freeness can be verified by an enumeration of the loose
tree and if the index tree is finite the freeness is independent of the
enumeration.  We do  not at first
claim so much for an infinite tree.  Rather we
introduce a notion of a locally free (loose) tree that has the
following properties.  If $\Mun$ is free under some enumeration then
$\Mun$ is locally free.  The property of being locally free does not
depend on the enumeration.  We will establish the existence of
compatibility prime
models over locally free loose trees (which have height at most
$\omega$).  From this we deduce that if such a tree is free is under
one enumeration then it is free under any enumeration.


\jtbdef  A loose tree of models $\Mun$ is {\em locally free in $N$} if
$\Mun=\langle M_t: t\in {\bf T}\rangle$
is contained in $N$ and for
every finite subtree
${\bf T}_1 \subseteq {\bf T}$  the finite loose tree
$\Mun|{\bf T}_1
=\langle M_t: t\in {\bf T_1}\rangle$ is free inside $N$.

\jtbnumpar{Remark}  This definition relies on our restriction to
subtrees of $\lambda^{<\omega}$.  Since subtrees are closed under
predecessor trees of greater height cannot be covered by finite
subtrees.  Thus, when we deal with trees of greater height we will
modify the definition of locally free but leave the meaning the same on
the low trees considered here.

The next proposition is obvious.
\begin{prop}
\label{freeenumimplieslocallyfree}
If the loose tree $\Mun$ is free inside for $N$ for some
enumeration then $\Mun$ is locally free inside $N$.
\end{prop}
Now we define `prime' models over locally free loose trees.  Definition
\ref{defprimloosetree}
is almost a special case of this.  (Ostensibly, a map could take a free
tree to a locally free tree so prime over locally free is more
restrictive than prime over free.)
\jtbdef
\label{defprimlocloosetree}
The model $M$ is {\it compatibility prime over the
locally free
loose tree
$\Mun$} via the
embedding $f$ of $\Mun$ into $M$
if the image of
$\Mun$ under $f$ is locally free inside $M$
and for every $N \in K$
and
$g$ embedding $\Mun$ into $N$
with the image of $g$ locally free inside $N$
such that $M$ and $N$ are compatible over $\Mun$ via $f$ and $g$
there is a $\hat g$ mapping $M$ into $N$
such that $g = \hat g \circ f$.

The name of this notion is unconscionably long so we resort to the
following abbreviation:
$\LFP(\Mun,M,f)$.  We omit the $f$ if $\Mun$ is concretely realized
in $M$ (whence $f$ is family of inclusions.)
\endjtbdef
We are going to show that if a loose tree $\Mun$ is locally free there
is a `prime' model over the loose tree
$\Mun$.
We need the following notation (extending \ref{treenot} and
\ref{notformainlemma}) for trees obtained from a representation of
given tree as a union of (smaller) trees
to carefully state and  prove this
proposition.

\jtbnot
\label{notformainlemma1}
Suppose ${\bf T} = \union_{\alpha < \lambda} I_{\alpha}$ where
the $I_{\alpha}$ form an
increasing continuous sequence of ideals in
{\bf T}.
\begin{enumerate}
\item
$\Mun_{\alpha}$
denotes the quotient
of $\Mun$ indexed by the quotient tree
$\{t:t\in
{\bf T}_{I_{\alpha}}\}$.
\item
For $\beta < \alpha$,
$\Mun_{\alpha,\beta}$
denotes the quotient of $\Mun|(I_{\alpha}$ by $I_{\beta}$.
That is,
$\Mun_{\alpha,\beta}$  is
indexed by the quotient tree $\{t:t\in {(I_{\alpha}})_{I_{\beta}}\}$
(Remember $I_{\alpha}$ is a subtree of ${\bf T} $ so this is
a natural extension of our
previous notation \ref{notformainlemma}.)
\end{enumerate}

Thus $\Mun_{\alpha}$ just abbreviates $\Mun_{I_{\alpha}}$.
Although $\Mun_{\alpha}$ is a quotient of $\Mun$ and
$\Mun_{\alpha,\beta}$ is a quotient of $\Mun|(I_{\alpha})$, as sets
of models $\Mun_{\alpha}$ and
$\Mun_{\alpha,\beta}$ are contained in $\Mun$.
\endjtbnot

Here is the main result of this section.
Condition 1 is the result we really want
(see Conclusion~\ref{conclusion})
but to establish it we must prove Conditions 1 and 2 by a simultaneous
induction.  A first try to prove this theorem by induction would
decompose a tree
{\bf T} as a union of properly smaller ideals and then choose a
`prime' model over each of these ideals and take the canonically
prime model over this chain of models.  The difficulty is to
guarantee that the sequence of models forms a chain.  The double
induction accomplishes this end.
\begin{thm}
\label{primoverdiag}
Suppose
$\Mun$ is a locally free
loose tree indexed by ${\bf T_0}$ in $N$.
Conditions
1) and 2) will be
proved by simultaneous
induction on $\lambda$.
\begin{description}
\item[Condition 1.]
If $|{\bf T}_0| = \lambda$ then
there is a compatibility prime model for locally free
loose trees over $\Mun$.
That is, there is an $N$ satisfying $\LFP(\Mun,N)$.
\item[Condition 2.]
If $I$ is an ideal of ${\bf T_0}$ and  $|I| = \lambda$
then
$\Mun_I(M_I/\emptynode)$ is
locally free inside $M$.
\end{description}
\end{thm}

Note that the cardinality of ${\bf T} _0$ is not bounded in condition
2.  The following lemma provides most of the technicalities of the
proof.
\begin{lemm}[The key lemma]
\label{techlemmalocfree}
Suppose
$\Mun$ is a free
loose tree in $N'$
indexed by the finite tree
{\bf T} and suppose $s$ is a minimal ($\neq \emptynode$)
element of {\bf T}.  Suppose also
$\langle M^{\alpha}_s:\alpha \leq \delta \rangle$ is an increasing
{\bf K}-continuous sequence satisfying the following conditions.
\begin{enumerate}
\item
For each $\alpha < \delta$,
$\Mun^{\alpha} = \Mun(M^{\alpha}_s /s)$ is
a free loose tree inside $N'$.
\item
For each $\alpha < \delta$,
$M_{\emptynode} \subm M^{\alpha}_s$.
\item
If $t^- = s$ then $\ind
{M_t} {{M_t} \inter M^{0}_s} {M^{\alpha}_s}$.
\end{enumerate}
Then
$\Mun(M^{\delta}_s/s)$
is free inside $N'$.
\end{lemm}

\PROOF.  Fix an enumeration $\tbar$ of {\bf T} such that the elements
not in the cone with vertex $s$ come first.  Suppose $s = t_k$ and
$|{\bf T} | = n$. Let $\Nun^{\alpha}$
witness that $\Mun^{\alpha}$ is free.  By Corollary~\ref{rutcor} we may
assume that $N^{\alpha}_i
= N^{0}_i$ if $i < k$.  (This use of Corollary~\ref{rutcor} only
simplifies notation; the later use is essential.)
In order to discuss uniformly the trees $\Mun^{\alpha}$, we refer below
to models $M^{\alpha}_t$.  Unless $t = s$, $M^{\alpha}_t = M^0_t = M_t$.

Expand the tree {\bf T} to $\hat {\bf T} $ by adding a new element $r$
with $r^-=s$ but with no elements above $r$.  We define for
each $\alpha\leq\delta$ a loose tree of models
${\hat {\Mun}}^{\alpha}$ and a sequence of models ${\hat {\Nun}}^{\alpha}$.
 \[ \hat M ^{\alpha}_x  = \left \{ \begin{array}{ll}
                    M^0_{x} & \mbox{if $x \in {\bf T}$} \\
                    M^{\alpha}_{s} & \mbox{if $x=r$}
                   \end{array}
          \right. \]
\[ \hat N ^{\alpha}_i  = \left \{ \begin{array}{ll}
                    N^0_i & \mbox{if $i \leq k$} \\
                    N^{\alpha}_{i-1} & \mbox{if $i > k$}
                   \end{array}
          \right. \]

Consider the following enumeration $\ubar$ of $\hat {\bf T} $ which we
obtain from $\tbar$.
\[ u_i  = \left \{ \begin{array}{lll}
                    t_i & \mbox{if $i \leq k$} \\
                    r   & \mbox{if $i = k+1$} \\
                    {t_{i-1}} & \mbox{if $k+1 < i \leq n-1$}
                   \end{array}
          \right. \]
Recall that if $i \neq k$,
$M^0_{t_i}
=M^{\alpha}_{t_i}$.  Note that both $M^0_s$ and $M^{\alpha}_s$
occur in ${\hat {\Mun}} ^{\alpha}$
($\hat M ^{\alpha}_s = M^0_s$,
$\hat M^{\alpha}_r =
\hat M^{\alpha}_{u_{k+1}} =
M^{\alpha}_s$).
We will show below that for each $\alpha$, ${\hat {\Nun}} ^{\alpha}$
witnesses that ${\hat {\Mun}} ^{\alpha}$ is free inside $N'$.
Assuming this
fact we now complete the proof of the lemma.

Since we proved that the freeness of a loose tree does not depend on the
enumeration, each
${\hat {\Mun}} ^{\alpha}$
is also free by an
enumeration
that places $r$ last.  That is
${\hat {\Mun}} ^{\alpha}$ is free
with respect to the eumeration $\vbar$
defined as follows.
\[v_i  = \left \{ \begin{array}{lll}
                    t_i & \mbox{if $i \leq n-1$} \\
                    r   & \mbox{if $i = n$}
                   \end{array}
          \right. \]
Let $I = \{v_i:i < n\}$. By Corollary~\ref{rutcor} any sequence, in
particular ${\hat {\Nun}} ^{0}= \langle N^0_0 \ldots N^0_{n-1}\rangle$,
which
witnesses the freeness of ${\hat {\Mun}} ^{\alpha}|I$ can be extended to a
sequence witnessing the freeness of ${\hat {\Mun}} ^{\alpha}$.  Let $N$
denote $N^0_{n-1}$.  (It was to make the choice of $N$ independent of
$\alpha$ that we needed to prove Theorem~\ref{rut} and
Corollary~\ref{rutcor}.) Now for each $\alpha$
since ${\hat {\Mun}} ^{\alpha}$ is free with respect to $\vbar$
we have
$$\ind N {\hat M^{\alpha}_{r} \inter \, \hat M^{\alpha}_{r^-}}
{\hat M^{\alpha}_{r}}$$
inside $\hat N^{\alpha}_n$.
That is,
$$\ind N {M^{\alpha}_{s} \inter M^0_{s}}
{M^{\alpha}_{s}}$$
inside some $N^{\alpha}_n$.
This implies
and since $M^0_{s} \subseteq
{M^{\alpha}_{s}}$ that for each $\alpha$
$$\ind N { M^0_{s}}
{M^{\alpha}_{s}}$$
inside $N'$.
Therefore, by L1,
$$\ind N { M^0_{s}}
{M^{\delta}_{s}}$$
inside $N'$.  Choosing $N_n$ prime over $N \union
{M^{\delta}_{s}}$, we finish the proof of the key lemma from the
assumption.

We are left with verifying the assumption; that is,
showing that for each $\alpha$, ${\hat {\Nun}}
^{\alpha}$ witnesses that ${\hat {\Mun}} ^{\alpha}$ is free inside $N'$
with respect to $\ubar$.
We will rely on the following fact.  For each $\alpha$
and each $i < n$
\begin{equation}
\label{techeqn1}
\ind
{M^{\alpha}_{t_i}}
{M^{\alpha}_{t_i} \inter M^{\alpha}_{t^-_i}}
{N^{\alpha}_{i-1}}
\hbox{\rm \/
inside $N^{\alpha}_i$}.
\end{equation}

We must show that for each $i < n+1$,
\begin{equation}
\label{techeqn2}
\ind
{\hat M^{\alpha}_{u_i}}
{\hat M^{\alpha}_{u_i} \inter \hat M^{\alpha}_{u^-_i}}
{\hat N^{\alpha}_{i-1}}
\hbox{\rm \/
inside $\hat N^{\alpha}_i$}.
\end{equation}
For $i \leq k$ this is an immediate translation of property
\ref{techeqn1} (with $\alpha = 0$).
For $i = k+1$ property \ref{techeqn2} translates to
\begin{equation}
\label{techeqn3}
\ind
{M^{\alpha}_{s}}
{M^{\alpha}_{s} \inter  M^{0}_{s}}
{N^{0}_{k}}
\hbox{\rm \ inside $N^{\alpha}_k$}.
\end{equation}
since
$\hat M^{\alpha}_{u^-_{k+1}} =
\hat M^{\alpha}_{u_{k}} =
M^{0}_{t_{k}}=
M^{0}_{s}$.
Remember that
$N^0_k$ is prime over $N^0_{k-1} \union M^0_s$
and since
$N^0_{k-1}=
N^{\alpha}_{k-1}$,
$N^{\alpha}_k$ is prime over $N^0_{k-1} \union M^{\alpha}_s$.
Moreover,
$M^{\alpha}_{s}
\inter
M^{\alpha}_{s-} =
M^{\alpha}_{s}
\inter
M^{\alpha}_{\emptynode}=
M^{\alpha}_{\emptynode}$.
Now property \ref{techeqn3} follows by base extension
from
\begin{equation}
\label{techeqn4}
\ind
{M^{\alpha}_{s}}
{M^{\alpha}_{s} \inter  M^{\alpha}_{\emptynode}}
{N^{\alpha}_{k-1}}
\hbox{\rm \/
inside $N^{\alpha}_k$}
\end{equation}
which is an instance of property \ref{techeqn1}.

For $i = k+2$, since $u^-_{k+2} = t^-_{k+1} = t_k = s$
property \ref{techeqn2} becomes
\begin{equation}
\label{techeqn5}
\ind
{M^{\alpha}_{t_{k+1}}}
{M^{\alpha}_{t_{k+1}} \inter  M^{0}_{s}}
{N^{\alpha}_{k}}
\hbox{\rm \/
inside $N^{\alpha}_{k+1}$}.
\end{equation}
Now by assumption iii) of the Lemma,
$\ind {M^{\alpha}_{t_{k+1}}} {{M^{\alpha}_{t_{k+1}}} \inter M^{0}_s}
{M^{\alpha}_s}$
so we can deduce property \ref{techeqn5} from property \ref{techeqn1}
and transitivity of independence (Theorem~\ref{transind}).

For $i > k+2$, property \ref{techeqn2} translates to
\begin{equation}
\label{techeqn6}
\ind
{M^{\alpha}_{t_{i-1}}}
{M^{\alpha}_{t_{i-1}} \inter \, M^{\alpha}_{t^-_{i-1}}}
{N^{\alpha}_{i-2}}
\hbox{\rm \/
inside $ N^{\alpha}_i$}.
\end{equation}
which follows immediately from property \ref{techeqn1}.

We will apply the following consequence of this lemma
(with $I_{\delta}$  a proper ideal of {\bf T} )
in the main proof.

\begin{cor}
\label{cortotechlemmlocfree}
Suppose $\Mun$ is a locally free loose tree  in some $N$ which is
indexed by {\bf T}
and $\langle
I_{\alpha}:\alpha \leq \delta\rangle$ is a continuous increasing
sequence of ideals in {\bf T}.  Suppose further that there is a {\bf
K}-continuous sequence $\langle N_{\alpha}:\alpha \leq \delta \rangle$
such that for each $\alpha <
\delta$, $\union_{t \in I_{\alpha}}M_t \subseteq N_{\alpha} \subm N$ and
$\Mun_{\alpha}(N_{\alpha}/\emptynode)$
is a locally free loose tree
inside $N$.
Then
$\Mun_{\delta}(N_{\delta}/\emptynode)$ is a locally free loose tree
inside $N$.
\end{cor}
\PROOF. Without loss of generality, (since we are trying to establish
local freeness),
${\bf T} -I_{\delta}$
is finite.  Let $\{t_i:i < k\}$ be the
minimal elements of
${\bf T} -I_{\delta}$
and let $s_i \in I_{\delta}$ be the predecessor of $t_i$ in
$I_{\delta}$.
Again, without loss of generality, we may assume all the $s_i$ are in
$I_0$.
For each $\alpha \leq \delta$, we define a loose tree
$\tilde {\Mun}^{\alpha}$ (not following our previous conventions)
as follows.
$\tilde {\Mun}^{\alpha}$
is indexed by
$({\bf T} -I_{\delta})
\union \{s,\emptynode\}$ where $s$ is
interpolated above $\emptynode$ and below each member of
${\bf T} -I_{\delta}$.
$\tilde M^{\alpha}_s = M_t$ if
$t \in ({\bf T} -I_{\delta})$,
$N_0$ if $t =\emptynode$ and $N_{\alpha}$ if $t =s$.
$\tilde {\Mun}^{\alpha}$
is a
free loose tree
inside $N$.  (Just enumerate $\emptynode$, $s$ and then the rest using
$\Mun_{\alpha}(N_{\alpha}/\emptynode)$
is a locally free loose tree
inside $N$.)
By the choice of  the $t_i$ and $s_i$, and the local freeness of
$\Mun_{\alpha}(N_{\alpha}/\emptynode)$ inside $N$
we verify the third assumption of
Lemma~\ref{techlemmalocfree}.  In the present context
(since $M_{s_i} \subseteq N_0$ and $M_{t_i} \inter M_{s_i}\subm N_0$)
it translates as
for each $i < k$,
$\ind
{M_{t_i}} {{M_{t_i}} \inter N_0} {N_{\alpha}}$.
Now by  Lemma~
\ref{techlemmalocfree}
$\tilde {\Mun}^{\delta}$
is free.  Applying the omission lemma to
$\tilde {\Mun}^{\delta}$
we conclude that
$\Mun_{\delta}(N_{\delta}/\emptynode)$ is a locally free loose tree
inside $N$.
\endproof

Now we continue with the proof of the main theorem
(Theorem
\ref{primoverdiag}).
\jtbnumpar{Proof of
\ref{primoverdiag}}
We assume by induction that Condition 1 of the theorem holds for
any tree with cardinality less than $\lambda$ and
Condition 2 of the theorem holds for any tree and any ideal of
cardinality less than $\lambda$.

Fix a subtree ${\bf T}_1 \subseteq {\bf T} _0$
with $|{\bf T} _1|= \lambda$.
We will describe
a construction relative to
${\bf T} _1$.
Then to verify each of
conditions 1 and 2 we use a different choice of ${\bf T}_1$.

 As in \ref{notformainlemma1},
suppose ${\bf T_1} = \union_{\alpha < \lambda} I_{\alpha}$ where
the $I_{\alpha}$ are an increasing continuous sequence of ideals in
${\bf T}_1$ and $|I_{\alpha}| <\lambda$.
We define by induction on $\alpha < \lambda$, an increasing
{\bf K}-continuous sequence $\langle N_{\alpha}:\alpha <\lambda\rangle$ of
{\bf K}-submodels of $N$ satisfying the following conditions.
\begin{enumerate}
\item
If $\alpha = \beta + 1$ then $N_{\alpha}$ is
prime for loose trees
over
$\Mun_{\alpha,\beta}(N_{\beta}/\emptynode)$.
\item
$\Mun_{\alpha}
(N_{\alpha}/\emptynode)$ is locally free.
\end{enumerate}
Recall that
$\Mun_{\alpha,\beta}$ denotes the quotient of $\Mun|(I_{\alpha})$ by
$I_{\beta}$.
Note that
$\Mun_{\alpha}$ is indexed by
$({\bf T}_0 - I_{\alpha})
\union \{\emptynode\}$
(not
$({\bf T}_1 - I_{\alpha})\union \{\emptynode\}$
).  This is crucial for the verification of
Condition 2.

There are three cases in the construction.
\begin{description}
\item[
$\alpha = 0$.]
Condition i) doesn't apply.  Condition ii) holds
by the induction hypothesis applied to Condition 2) of the main
theorem.
\item[
$\alpha$ is a limit ordinal.]
Let $N_{\alpha}$ be canonically
prime over $\langle N_{\beta}:\beta <\alpha\rangle$.
Condition i) does not
apply and condition ii) is immediate from the last result: Corollary
\ref{cortotechlemmlocfree}.
\item[
$\alpha = \gamma+1$.]
Applying condition ii) to $\gamma$,
$\Mun_{I_{\gamma}}(N_{\gamma}/\emptynode)$ is locally free in $N$.
So the
subdiagram
$\Mun_{\alpha,\gamma}$ is locally free
in $N$.  So by induction applied to
Condition 1) of the main theorem with $I_{\alpha}$ as $I$ there is
an $N_{\alpha}$ satisfying condition
i).  That is, $N_{\alpha}$  is
prime for loose trees
over
$\Mun_{\alpha,\gamma}(N_{\gamma}/\emptynode)$.  In particular, note
$N_{\gamma}$ is embedded into $N_{\alpha}$.
Applying the induction hypothesis for Condition 2 (with ${\bf T} _0 =
\Mun_{I_{\alpha}}$)
$\Mun_{\alpha}
(N_{\alpha}/\emptynode)$ is locally free
and we satisfy condition ii).
\end{description}
This completes the construction.

To see that Condition 1 holds of ${\bf T} _0$, take ${\bf T}_1$ as
${\bf T}_0$.
We must show that if $N_{\lambda}$ is taken
canonically
prime over $\langle N_{\beta}:\beta <\lambda\rangle$
then $N_{\lambda}$ satisfies $\LFP(\Nun,N_{\lambda})$.
Suppose $f= \langle f_t:t \in {\bf T} \rangle$ maps
the loose tree $\Mun$
isomorphically to a loose tree $\Mun'$ which is locally free in some
$M'$.  We must extend $f$ to an embedding of $N_{\lambda}$ into $M'$.
\par
For this, repeat the preceding argument constructing a sequence
$N'_{\alpha}$ for $\alpha < \lambda$ inside $M'$ which satisfy
conditions i) and ii) (for $\Mun'$) and simultaneously construct
a continuous increasing sequence of maps
$g_{\alpha}$ for $\alpha < \delta$ which map $N_{\alpha}$ isomorphically
onto $N'_{\alpha}$ and such that $g_\alpha$ extends each $f_t$ with
$t\in I_{\alpha}$.
We conclude Condition 1.
\par

For Condition 2
we must show that
if $I$ is an ideal of {\bf T} and  $|I| = \lambda$ then $\Mun_I$ is
locally free.
Applying the construction above with $I$ as ${\bf T}_1$,
we have constructed the $N_{\alpha}$ so that
$\Mun_{\alpha}(N_{\alpha}/\emptynode)$ is locally free.  Thus, by
Corollary \ref{cortotechlemmlocfree}
$\Mun_{\alpha}(N_{\delta}/\emptynode)$ is locally free and we can
conclude Condition 2.
\endproof
\jtbnumpar{Conclusion}
\label{conclusion}
Let {\bf K} be an adequate class and suppose
$\Mun$ is a loose tree of {\bf K} -models indexed by a subtree of
$\lambda^{<\omega}$.
If $\Mun$ is locally free then there is a compatibility prime model $M$
over $\Mun$ (i.e.
$\LFP(\Mun,M,f)$, see \ref{defprimlocloosetree}.)

In fact, if $\Mun$ (indexed by a subtree of $\lambda^{<\omega}$) is free
under one enumeration then it is free under any enumeration.  For, by
Proposition~\ref{freeenumimplieslocallyfree}
it is locally free.  Let $\tbar=\langle t_i;i<|{\bf T}|\rangle$ be an
arbitrary enumeration of $\Mun$.  For each $\gamma < |{\bf T}|$ let
$I_{\gamma}$ be the ideal containing each $t_i$ with $i < \gamma$.
Now construct a family of models $\langle N_{\gamma}:\gamma < |{\bf T}|
\rangle$ as in the proof of Theorem~\ref{primoverdiag}.  By the
basic properties of independence (using L1 at limit stages) these
models witness that $\Mun$ is free under the given enumeration.

We have shown  that if $\Mun$ is indexed by a subtree of
$\lambda^{<\omega}$ and if $\Mun$ is free with respect to some
enumeration then there is a compatibility prime model over $\Mun$.  In
\cite{BaldwinShelahprimalii} we will show that if {\bf K} is at all
manageable then there is only one compatibility class over $\Mun$.
We are attempting to show that models in manageable classes are
`tree-decomposable' by trees of small height.  The existence of such
prime models is an essential step in this program.
\bibliography{ssgroups}
\bibliographystyle{plain}
 \end{document}